\documentclass[11pt,intlimits,reqno]{amsart}

\usepackage{hyperref}
\usepackage{graphicx}
\usepackage{color}
\usepackage{amsfonts,amssymb}
\usepackage{amsthm}

\usepackage{a4wide}

\newtheorem{neu}{}[section]
\newtheorem{Cor}[neu]{Corollary}
\newtheorem*{Cor*}{Corollary}
\newtheorem{Thm}[neu]{Theorem}
\newtheorem*{Thm*}{Theorem}
\theoremstyle{definition}
\newtheorem{Prop}[neu]{Proposition}
\newtheorem*{Prop*}{Proposition}
\newtheorem{Lemma}[neu]{Lemma}
\newtheorem*{Rmk*}{Remark}
\newtheorem{Rmk}[neu]{Remark}
\newtheorem{Ex}[neu]{Example}
\newtheorem*{Ex*}{Example}

\newtheorem{Def}[neu]{Definition}

\newcommand{\N}{\mathbb{N}}

\newcommand{\Z}{\mathbb{Z}}
\newcommand{\R}{\mathbb{R}}
\newcommand{\C}{\mathbb{C}}
\newcommand{\CP}{\C\mathrm{P}}

\newcommand{\zw}{\hspace*{0.75em}}
\newcommand{\pf}{\longrightarrow}
\newcommand{\Pf}{\Longrightarrow}

\newcommand{\wrt}{with respect to }

\newcommand{\CZ}{\mu_{\mathrm{CZ}}}
\newcommand{\Mas}{\mu_{\mathrm{Maslov}}}
\newcommand{\Morse}{\mu_{\mathrm{Morse}}}
\newcommand{\db}{\bar{\partial}_{J,H}\,}   

\newcommand{\im}{\mathrm{im\,}}

\newcommand{\ind}{\mathrm{ind\,}}

\newcommand{\om}{\omega}
\newcommand{\ev}{\mathrm{ev}}

\newcommand{\Poincare}{Poincar\'{e} }

\newcommand{\A}{\mathcal{A}}

\newcommand{\Pt}{\widetilde{\mathcal{P}}}

\newcommand{\F}{\mathcal{F}}
\newcommand{\U}{\mathcal{U}}

\newcommand{\M}{\mathcal{M}}
\newcommand{\Mh}{\widehat{\mathcal{M}}}

\renewcommand{\L}{\mathcal{L}}
\newcommand{\Lt}[1]{\widetilde{\mathcal{L}#1}}
\renewcommand{\H}{\mathrm{H}}

\newcommand{\PSS}{\mathrm{PSS}}
\newcommand{\Ham}{\mathrm{Ham}}
\newcommand{\Nov}{\Lambda(\Gamma)}
\newcommand{\Novo}{\Lambda_0(\Gamma)}
\newcommand{\PD}{\mathrm{PD}}
\newcommand{\CF}{\mathrm{CF}}
\newcommand{\HF}{\mathrm{HF}}
\newcommand{\CM}{\mathrm{CM}}
\newcommand{\HM}{\mathrm{HM}}
\newcommand{\QH}{\mathrm{QH}}
\newcommand{\Crit}{\mathrm{Crit}}

\newcommand{\beq}{\begin{equation}}
\newcommand{\beqn}{\begin{equation}\nonumber}
\newcommand{\eeq}{\end{equation}}

\newcommand{\bea}{\begin{equation}\begin{aligned}}
\newcommand{\bean}{\begin{equation}\begin{aligned}\nonumber}
\newcommand{\eea}{\end{aligned}\end{equation}}

\numberwithin{equation}{section}

\begin{document}
\title[On the extrinsic topology of Lagrangian submanifolds]{On the extrinsic topology of\\Lagrangian submanifolds}
\author{Peter Albers}
\address{
    Mathematisches Institut\\
    Universit\"at Leipzig\\
    04109 Leipzig\\
    Germany}
\email{albers@math.uni-leipzig.de}
\date{April 2005}
\begin{abstract}
We investigate the extrinsic topology of Lagrangian submanifolds and of their submanifolds in closed symplectic manifolds using Floer homological
methods. The first result asserts that the homology class of a displaceable monotone Lagrangian submanifold vanishes in the homology of the ambient
symplectic manifold.
Combining this with spectral invariants we provide a new mechanism for proving Lagrangian intersection results e.g.~entailing
that any two simply connected Lagrangian submanifold in $\CP^n\times\CP^n$ must intersect.
\end{abstract}
\maketitle
\tableofcontents

\section{Introduction}
The fact that a manifold $L$ admits an embedding into a symplectic manifold as a \textit{Lagrangian} submanifold yields restrictions on
the topology of $L$. This \textit{intrinsic topology} of Lagrangian submanifolds has been studied quite extensively, cf.~the
recent survey \cite{Biran_Geometry_of_symplectic_intersections}.

In this paper, we investigate the \emph{extrinsic topology} of Lagrangian submanifolds: We address
the question how a Lagrangian submanifold lies homologically in the ambient symplectic manifold. Among other things this turns into a new Lagrangian
intersection mechanism.

In the basic theorem \ref{thm:basic_theorem}, we provide for a monotone closed Lagrangian submanifold $L$
in a closed symplectic manifold $M$ a Floer-theoretic representation of the
homomorphisms $\iota_k:\H_k(L;\Z/2)\pf\H_k(M;\Z/2)$ for degrees $k>\dim L+1-N_L$, where $N_L$ denotes the
\emph{minimal Maslov number} of $L$ (see definition \ref{def:minimal_numbers}).
In the second basic theorem \ref{thm:basic_theorem_with_low_energy} we remove the restriction on the degree of $\iota_k $ for
Hamiltonian functions having sufficiently small Hofer norm.
If the Lagrangian submanifold $L$ is displaceable, theorem \ref{thm:basic_theorem} implies a
vanishing result for the homomorphisms $\iota_k:\H_k(L;\Z/2)\pf\H_k(M;\Z/2)$:
\begin{Thm*}
Let $(M,\om)$ be a monotone closed symplectic manifold and $L\subset M$ a monotone closed Lagrangian submanifold with $N_L\geq2$. If the Lagrangian
submanifold $L$ is (Hamiltonianly) displaceable, then the homomorphism $\iota_k$ vanishes for degrees $k>\dim L+1-N_L$.
\end{Thm*}

\begin{Cor*}
For any displaceable monotone Lagrangian submanifold $L$ with $N_L\geq2$,
\beqn
[L]=0\in\H_m(M;\Z/2)\,.
\eeq
\end{Cor*}

Applying the basic theorems \ref{thm:basic_theorem} and \ref{thm:basic_theorem_with_low_energy} to spectral capacities
leads to a new mechanism for producing Lagrangian intersection results. Namely, in a symplectic manifold with
finite spectral capacity all monotone Lagrangian submanifolds $L$ of minimal Maslov number $N_L>\dim L+1$
intersect each other, cf.~theorem \ref{thm:high_minimal_Maslov_number_implies_intersection}. An instance of this new mechanism is:

\begin{Cor*}
Any two simply connected Lagrangian submanifolds in $\CP^n\times\CP^n$ intersect.
\end{Cor*}

An example of a simply connected Lagrangian submanifold of $\CP^n\times\CP^n$
is the anti-diagonal $\bar{\Delta}=\{(\bar{z},z)\mid z\in\CP^n\}$.\\[1ex]
\textit{Organization of the paper.}\quad
In section \ref{sec:Floer_homology_and_spectral_capacities} we recall the construction of Floer homology for semi-positive symplectic manifolds and the
Piunikhin-Salamon-Schwarz isomorphism. Furthermore, we give the definition of the action selector and of spectral
capacities via the action filtration on Floer homology. Theorems \ref{thm:basic_theorem} and
\ref{thm:basic_theorem_with_low_energy} are stated in section \ref{sec:basic_theorems} and proved in section \ref{sec:Proof_of_the_basic_theorems}.
In section \ref{sec:applications} we derive various applications,
which are divided into two subsections, on the extrinsic topology and on Lagrangian intersections.\\[1ex]
\textbf{Acknowledgements.}\quad
The results of this paper are by-products of my Ph.D. thesis \cite{Albers_PhD} and are partly contained therein.
I would like to express my gratitude to my advisor Matthias Schwarz.

I was financially supported by the DFG through the Graduiertenkolleg
"Analysis, Geometrie und ihre Verbindung zu den Naturwissenschaften" at the University of Leipzig
and the Schwerpunktprogramm "Globale Differentialgeometrie", grant SCHW 892/2-1.
\section{Floer homology and spectral capacities}
\label{sec:Floer_homology_and_spectral_capacities}

Despite the fact that we consider here Floer homology only for monotone symplectic manifolds we recall briefly and without proofs the construction
of Floer homology for semi-positive symplectic manifolds with the help of Novikov rings. Even though Novikov rings can be avoided when dealing
with monotone symplectic manifold, the statements are much clearer when using Novikov rings. All details can be found in
\cite{
McDuff_Salamon_J_holomorphic_curves_and_symplectic_topology,
Salamon_lectures_on_floer_homology,
Hofer_Salamon_Floer_homology_and_Novikov_rings}.
\begin{Def}
A closed symplectic manifold $(M^{2m},\om)$ is called \textit{semi-positive} if
\beq\tag{SP}\label{property:semi_positive}
\om(A)>0,\;c_1(A)\geq 3-m\quad\Pf\quad c_1(A)\geq0
\eeq
for all $A\in\pi_2(M)$, where $c_1$ is the first Chern class of $(M,\om)$.
\end{Def}
Most symplectic manifolds treated here are \textit{monotone}, i.e.~they satisfy $\om|_{\pi_2(M)}=\lambda\cdot c_1|_{\pi_2(M)}$
for some constant $\lambda>0$. Monotone symplectic manifolds are clearly semi-positive.
\begin{Def}\label{def:minimal_numbers}
Let $(M,\om)$ be a symplectic manifold and $L\subset M$ a Lagrangian submanifold.
The \textit{minimal area of a non-constant holomorphic sphere} in $M$ is
\beqn
A_M=\inf\Big\{\;\int_{S^2}s^*\omega\;\big|\;s\text{ is a non-constant holomorphic sphere in }M\;\Big\}\,.
\eeq
The \textit{minimal area of a non-constant holomorphic disk} on $L$ is
\beqn
A_L=\inf\Big\{\;\int_{D^2}d^*\omega\;\big|\;d\text{ is a non-constant holomorphic disk with boundary on }L\;\Big\}\,.
\eeq\\[-2ex]
The infimum over the empty set is by definition $+\infty$.
Furthermore, the \textit{minimal Chern number} $N_M$ of $M$ is the positive generator of the image $c_1(\pi_2(M))\subset\Z$ of the first
Chern class $c_1$, and the \textit{minimal Maslov number} $N_L$ of $L$ is the positive generator of the image $\Mas(\pi_2(M,L))\subset\Z$
of the Maslov index $\Mas$.
In case that $c_1$ or $\Mas$ vanish we set $N_M=+\infty$ or $N_L=+\infty$.
\end{Def}
\begin{Rmk}\label{rmk:orientable_Lagrangian_minimal_Maslov_geq_2}
The Maslov index equals
$2c_1$ on spheres and is congruent modulo 2 to the first Stiefel-Whitney class $w_1(L)$ of $L$ evaluated on the boundary loop of a disk.
Since $w_1(L)=0$ if $L$ is orientable, the Maslov index is an even number for orientable Lagrangian submanifolds.
Due to the fact that we set $N_L=+\infty$ in case $\Mas\equiv0$ we conclude $N_L\geq2$ for all orientable Lagrangian submanifolds.
\end{Rmk}
\noindent\textbf{The Novikov ring.}\quad
Let $\Gamma$ be the finitely generated Abelian group
\beq
\Gamma:=\frac{\pi_2(M)}{\ker c_1\cap\ker\om}\,.
\eeq
A grading on $\Gamma$ is defined via $\deg:\Gamma\pf\Z$, $\deg(A):=2c_1(A)$.
We set $\Gamma_k:=\deg^{-1}(k)$ and define the \textit{Novikov ring} associated to $\Gamma$ by
\begin{gather}\nonumber\textstyle
\Lambda(\Gamma):=\bigoplus_k\Lambda_k(\Gamma)\,,\text{ where}\\
\Lambda_k(\Gamma):=\Big\{\sum_{A\in \Gamma_k}\lambda_Ae^A\;\big|
        \;\#\{A\in\Gamma_k\,|\,\lambda_A\not=0,\om(A)\leq c\}<\infty\;\forall c\in\R\Big\}\,,
\end{gather}
with coefficients $\lambda_A$ in $\Z/2$. The ring structure is given by
\bea
\Lambda_k(\Gamma)\times\Lambda_l(\Gamma)&\pf\Lambda_{k+l}(\Gamma)\,,\\[1ex]
\big(\sum_{A\in \Gamma_k}\lambda_Ae^A\big)\cdot\big(\sum_{B\in \Gamma_l}\mu_Be^B\big)
                &\mapsto\sum_A\sum_B(\lambda_A\cdot\mu_B)\,e^{A+B}\\[0.5ex]
                &\qquad=\sum_C\big(\sum_A\lambda_A\cdot\mu_{C-A}\big)\,e^{C}\;.
\eea
$\Novo$ is a field and, in particular, $\Nov$ is a vector space over $\Lambda_0(\Gamma)$.
\begin{Ex*}[symplectically aspherical case]
If the symplectic manifold $(M,\om)$ is symplectically aspherical, i.e. $\om|_{\pi_2(M)}=0$ and ${c_1}|_{\pi_2(M)}=0$, then
\beq
\Gamma=\Gamma_0=\{0\}\quad\text{and}\quad\Nov=\Lambda_0(\Gamma)=\Z/2.
\eeq
\end{Ex*}
\begin{Ex*}[monotone case]
If $(M,\om)$ is monotone, i.e.~$\exists\lambda>0$ such that $\om|_{\pi_2(M)}=\lambda\cdot c_1|_{\pi_2(M)}$, then $\Gamma_0=\{0\}$ and $\Gamma\cong\Z$.
This implies that
\beq
\Lambda_0(\Gamma)=\Z/2\quad\text{and}\quad\Nov\cong\big(\Z/2\big)[\!\!\:[q]\!\!\:][q^{-1}]\,,
\eeq
i.e.~$\Nov$ is isomorphic to the ring of Laurent series with coefficients in $\Z/2$.
\end{Ex*}
There are various types of Novikov rings for which Floer homology and quantum cohomology can be defined,
cf.~\cite[chapter 11.1]{McDuff_Salamon_J_holomorphic_curves_and_symplectic_topology}. We choose here the field $\Novo$.\\[1ex]
\noindent\textbf{Floer homology over $\Novo$.}\quad
The group $\Gamma$ gives rise to a covering $\Gamma\rightarrow\Lt{M}\rightarrow\L M$
of the space $\L M$ of \emph{contractible loops} in $M$. The elements of $\Lt{M}$
are represented by equivalence classes $\bar{x}\equiv[x,d_x]$, where $x\in\L M$ and $d_x$ is an extension of the
contractible loop $x$ to the unit disk $D^2$. Two pairs $(x,d_x)$ and $(y,d_y)$ are equivalent if
\beq
(x,d_x)\sim(y,d_y)\quad\Longleftrightarrow\quad x=y\,,\quad c_1(d_x\#\bar{d}_y)=0\,,\quad \om(d_x\#\bar{d}_y)=0\,,
\eeq
where $d_x\#\bar{d}_y$ is the sphere obtained by glueing the two disks along their common boundary $x=y$.
The group $\Gamma$ acts on $\Lt{M}$ by concatenating the disk $d_x$ with the sphere $A$
\beq
\bar{x}\#A\;\equiv\;[x,d_x]\#A\;:=\;[x,A\#d_x]\quad\text{for all }A\in\Gamma\text{ and }\bar{x}\in\Lt{M}\,.
\eeq
For a Hamiltonian function $H:S^1\times M\pf\R$ we define the \textit{action functional} $\A_H$ on $\Lt{M}$ by\\[-2ex]
\beq
\A_H([x,d_x]):=\int_{D^2}d_x^*\,\om-\int_0^1H\big(t,x(t)\big)\,dt
\eeq
and the \textit{Hamiltonian vector field} $X_H$ associated to $H$ by $\om(X_H,\cdot)=-dH$. The time-1-map $\Phi_H$ of the flow induced
by the Hamiltonian vector field $X_H$ is called a \textit{Hamiltonian diffeomorphism}.
The set of critical points of $\A_H$ is
\beq
\Pt(H)=\big\{[x,d_x]\in\Lt{M}\mid \dot{x}(t)=X_H\big(t,x(t)\big)\big\}
\eeq
and is graded by the Conley-Zehnder index $\CZ:\Pt(H)\pf\Z$, $[x,d_x]\mapsto\CZ([d_x])$. We abbreviate $\Pt_k(H):=\CZ^{-1}(k)$.
The covering $\Lt{M}$ is chosen in such a way that the action functional $\A_H$ and the Conley-Zehnder index $\CZ$
are $\R$ resp.~$\Z$-valued.

The action functional and the grading behave as follows under the action of the group $\Gamma$
\beq
\A_H(\bar{x}\#A)=\A_H(\bar{x})+\om(A)\quad\text{and}\quad\CZ(\bar{x}\#A)=\CZ(\bar{x})+\deg(A)\,.
\eeq
For a non-degenerate Hamiltonian function $H$, i.e.~if the graph of $\Phi_H$ is transverse to the diagonal $\Delta\subset M\times M$,
we define the \textit{Floer chain groups} by
\beq
\CF_k(H):=\Big\{\sum_{\bar{x}\in \Pt_k(H)}\!\!\!a_{\bar{x}}\bar{x}\;\big|
        \;\#\{\bar{x}\,|\,a_{\bar{x}}\not=0,\A_H(\bar{x})\geq c\}<\infty\;\forall c\in\R\Big\}\,.
\eeq
The Floer chain groups become finite-dimensional vector spaces over the field $\Novo$
if we set
\bea
\Big(\sum_{A\in \Gamma_0}\lambda_Ae^A\Big)\cdot\Big(\sum_{\bar{x}}a_{\bar{x}}\bar{x}\Big)
        &:=\sum_{\bar{x}}\sum_A\lambda_Aa_{\bar{x}}\;\bar{x}\#(-A)\\
        &\;=\sum_{\bar{x}}\big(\sum_A\lambda_Aa_{\bar{x}\#A}\big)\bar{x}\,,
\eea
(cf.~\cite{Hofer_Salamon_Floer_homology_and_Novikov_rings}). In particular, $\A_H\big(e^A\cdot\bar{x}\big)=\A_H(\bar{x}\#(-A))=\A_H(\bar{x})-\om(A)$,
i.e.~the Novikov action of $A$ decreases the value of the action functional by $\om(A)$.

The boundary operator is defined by counting \textit{connecting Floer cylinders}. We define
\beq\label{def:moduli_space_of_Floer_connecting_traj_M(x,y,;J,H)_semi-positive}
\M(\bar{x},\bar{y};J,H):=\left\{u\in C^\infty(\R\times S^1,M)\left|\;
    \begin{aligned}
        &\partial_su+J(u)\big(\partial_tu-X_H(t,u)\big)=0\\
        &u(-\infty)=x,\quad u(+\infty)=y\\
        &d_x\#u\sim d_y
    \end{aligned}
\right.\right\},
\eeq
where $d_x\#u\sim d_y$ is in the Novikov sense, namely $c_1(d_x\#u\#\bar{d}_y)=0=\om(d_x\#u\#\bar{d}_y)$,
and where $J$ is a \textit{compatible almost complex structure} on $M$.

\begin{Thm*}[Floer-Hofer-Salamon \cite{Floer_Hofer_Salamon_Transversality_in_elliptic_Morse_theory_for_the_symplectic_action}]
For generic choices of $H$ and $J$ the moduli space $\M(\bar{x},\bar{y};J,H)$ is a
finite-dimensional manifold of dimension
\beq
\dim\M(\bar{x},\bar{y};J,H)=\CZ(\bar{y})-\CZ(\bar{x})
\eeq
admitting a free $\R$-action if $x\not=y$. Moreover, the moduli space
\beq\label{eqn:moduli_space_Floer_cylinders}
\Mh(\bar{x},\bar{y};J,H):=\M(\bar{x},\bar{y};J,H)\big/\R
\eeq
is a finite set for relative index one, i.e.~if $\CZ(\bar{y})-\CZ(\bar{x})=1$,
and it can be compactified by adding broken solutions for relative index two.
\end{Thm*}

For relative index one we denote the (mod 2) number of elements by
\beq
n(\bar{x},\bar{y}):=\#_2\big(\Mh(\bar{x},\bar{y};J,H)\big)
\eeq
and define the \textit{boundary operator of the Floer complex} on generators by
\bea
\partial:\CF_k(H)&\pf\CF_{k-1}(H)\\[0.5ex]
\bar{y}\;\;&\mapsto\;\partial\bar{y}\;:=\!\!\!\!\sum_{\CZ(\bar{x})=\CZ(\bar{y})-1}\!\!n(\bar{x},\bar{y})\cdot\bar{x}\,.
\eea
We extend it linearly to a $\Z/2$--vector space homomorphism. The Floer equation yields
\beq
n(\bar{x},\bar{y})\not=0\quad\Longrightarrow\quad \A_H(\bar{x})\leq\A_H(\bar{y})\,,
\eeq
i.e.~our conventions imply that the value of the action functional and the Conley-Zehnder index increase along Floer cylinders.
Furthermore, the following compactness result holds
\beq
\sum_{c_1(A)=0,\,\om(A)\leq c}\!\!\!\!\!\!\#\Mh(\bar{x},\bar{y}\#A;J,H)<\infty
\eeq
for all $\bar{x},\bar{y}\in\Pt(H)$ with $\CZ(\bar{y})-\CZ(\bar{x})=1$ and every $c\in\R$
(cf.~\cite{Hofer_Salamon_Floer_homology_and_Novikov_rings,
McDuff_Salamon_J_holomorphic_curves_and_symplectic_topology}).
These two facts imply that $\partial\bar{y}$ is a well-defined element in $\CF_{k-1}(H)$. The identity $n(\bar{x}\#A,\bar{y}\#A)=n(\bar{x},\bar{y})$
implies that $\partial$ actually is a $\Novo$-linear homomorphism. Floer's fundamental theorem asserts $\partial\circ\partial=0$,
so that \textit{Floer homology groups} are defined and form $\Novo$ vector spaces
\beq
\HF_*(H):=\mathrm{H}_*(\CF_*(H),\partial)\,.
\eeq
\begin{Ex*}[monotone case]
We recall that if $(M,\om)$ is monotone, then
\beq
\Gamma_0=\{0\}\quad\text{and}\quad\Lambda_0(\Gamma)=\Z/2\,.
\eeq
In particular, $\Pt_k(H)$ is finite, although in general $\Pt(H)$ will be infinite.
Furthermore, each Floer chain group $\CF_k(H)$ is a finite-dimensional vector space over the field $\Novo=\Z/2$.
In this case the Novikov conditions are empty. Nonetheless the full Novikov
ring $\Nov\cong\big(\Z/2\big)[\!\!\:[q]\!\!\:][q^{-1}]$ appears for instance in the following theorem by Hofer-Salamon.

For the moment let us denote the Floer homology over the full Novikov ring $\Nov$
by $\overline{\HF}_*(H)$, which is not $\Z$ but $\Z/(2N_M)$-graded, where $N_M$ is the minimal Chern
number of $M$ (cf.~definition \ref{def:minimal_numbers}). Both homology groups are related by the isomorphism
\beq
\overline{\HF}_{\bar{k}}(H)\;\cong\;\HF_k(H)\otimes_{\Novo}\Nov \quad\forall\;\bar{k}\in\Z/(2N_M)\,.
\eeq

\begin{Thm*}[Hofer-Salamon \cite{Hofer_Salamon_Floer_homology_and_Novikov_rings}]
There exists a $\Nov$-module isomorphism
\beq
\overline{\HF}_{\bar{k}}(H)\;\cong_{\Nov}
    \!\!\bigoplus_{i\equiv k\;\mathrm{mod}\;2N_M}\!\!\!\!\H^{m-i}(M;\Z/2)\otimes_{\Z/2}\Novo\,.
\eeq
\end{Thm*}
\end{Ex*}
\noindent\textbf{The Piunikhin-Salamon-Schwarz isomorphism.}\quad
In this section, we recall the construction of the \textit{Piunikhin-Salamon-Schwarz isomorphism},
which we call $\PSS$-isomorphism for brevity. The construction is presented
in \cite{Piunikhin_Salamon_Schwarz_Symplectic_Floer_Donaldson_theory_and_quantum_cohomology}
for semi-positive symplectic manifolds.

Since the $\PSS$-isomorphism plays a crucial role in the basic theorems below, we give here a fairly detailed exposition.

The $\PSS$-isomorphism is defined via counting solutions of a special boundary value problem whose solutions we call
\textit{plumber's helper solutions}. The corresponding moduli space $\M^{\mathrm{PSS}}(q,\bar{x};J,H,f,g)$
consists of pairs $(\gamma,u)$ of maps
\begin{gather}
\gamma:(-\infty,0]\pf M\zw\text{and}\zw u:\R\times S^1\pf M\zw \text{with }E(u)=\!\!\int_{-\infty}^\infty\int_0^1|\partial_su|^2\,dt\,ds<\infty\\[1ex]
\text{solving}\qquad\dot{\gamma}+\nabla^{g} f\circ\gamma=0
        \quad\text{and}\quad\partial_su+J(u)\big(\partial_tu-\beta(s)X_{H}(t,u)\big)=0\,,
\end{gather}\\[-1ex]
where $\beta:\R\rightarrow[0,1]$ is a smooth cut-off function satisfying $\beta(s)=0$ for $s\leq\tfrac12$ and $\beta(s)=1$ for $s\geq1$. The function
$f:M\rightarrow\R$ is a Morse function on $M$ and $\nabla^{g}$ is the gradient
\wrt a Riemannian metric $g$. Again $J$ is a compatible almost complex structure.
The pair $(\gamma,u)$ is required to satisfy boundary conditions
\beq
\gamma(-\infty)=q\,,\quad\gamma(0)=u(-\infty)\quad\text{and}\quad u(+\infty)=x\,,
\eeq
(see remark \ref{rmk:cut-off_implies_extension_u(infty)} below) where $q\in \Crit(f)$ is a critical point of $f$ and
$\bar{x}=[x,d_x]\in\Pt(H)$. Furthermore, we require the homotopy condition $\om(u\#\bar{d}_x)=0=c_1(u\#\bar{d}_x)$.
For brevity we denote this moduli space by $\M^{\mathrm{PSS}}(q,\bar{x})$. A pair $(\gamma,u)$ forms a plumber's helper.
\begin{figure}[h]
\begin{center}
\input{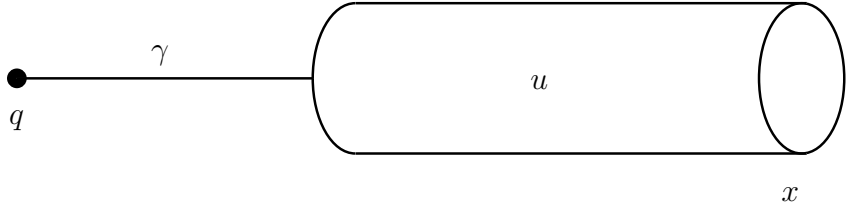}
\caption{A plumber's helper solution $(\gamma,u)\in\M^{\mathrm{PSS}}(q,\bar{x})$.}\label{fig_plumbers_helper}
\end{center}
\end{figure}
\begin{Rmk}\label{rmk:cut-off_implies_extension_u(infty)}
$\M^{\mathrm{PSS}}(q,\bar{x})$ is composed of gradient flow half-lines $\gamma$ for $f$ and Floer cylinders $u$.
Because of the cut-off function $\beta$ the cylinder $u$ actually is holomorphic for $s\leq\tfrac12$.
Since we impose the finite energy condition $E(u)<\infty$,
the punctured holomorphic disk $u\big|_{(-\infty,1/2)\times S^1}$ has
a removable singularity at the origin (cf.~\cite{McDuff_Salamon_J_holomorphic_curves_and_symplectic_topology}),
i.e.~$u$ has a continuous extension $u(-\infty)$. Therefore, the second boundary condition and the homotopy condition is meaningful.
\end{Rmk}

For generic choices of $H,J,f$ and $g$ the moduli space $\M^{\mathrm{PSS}}(q,\bar{x})$ is a smooth
manifold and
\beq
\dim \M^{\mathrm{PSS}}(q,\bar{x})= \CZ(\bar{x})+\Morse(q)-m\,,
\eeq
according to our normalization that for an autonomous $C^2$-small Morse function $f$ we have
$\CZ(d_x)=m-\Morse(x)$ for $\bar{x}=[x,d_x\equiv x]\in\Pt(f)$, i.e.~$x\in\Crit(f)$. As always $\dim M=2m$.

A standard computation for an element $(\gamma,u)\in\M^{\mathrm{PSS}}(q,\bar{x})$ shows\\[-2ex]
\beq\label{eqn:energy_of_plumbers_helper_solutions}
0\leq E(u)\leq\A_{H}(\bar{x})+\int_0^1\sup_M H(t,\cdot)dt\,.
\eeq
This implies a universal energy bound on plumber's helper solutions. In particular, sequences converge
up to breaking and bubbling.

Assumption \eqref{property:semi_positive} ensures compactness of
the 0-dimensional components of $\M^{\mathrm{PSS}}(q,\bar{x})$.
We denote by $\CM^{m-*}(f)$ the Morse co-chain complex associated to $f$ (and $g$), and define
\bea
\PSS:\CF_*(H)&\pf\CM^{m-*}(f)\\
\bar{x}\;\;&\mapsto\!\!\!\!\sum_{\Morse(q)=m-\CZ(\bar{x})}\#_2\M^{\mathrm{PSS}}(q,\bar{x})\cdot q\,.
\eea
\begin{Thm}[\cite{Piunikhin_Salamon_Schwarz_Symplectic_Floer_Donaldson_theory_and_quantum_cohomology}]
The 1-dimensional components of the moduli space $\M^{\mathrm{PSS}}(q,\bar{x})$ can be compactified by adding either triples $(\gamma,u,u')$, where
$(u,u')$ are broken Floer cylinders, or $(\gamma',\gamma,u)$, where $(\gamma',\gamma)$ are broken Morse trajectories.

This implies that $\PSS$ commutes with the Morse respectively the Floer (co)-boundary operator, i.e.
$\PSS$ descends to homology
\beq
\PSS:\HF_*(H)\pf\HM^{m-*}(f)\,.
\eeq
\end{Thm}
For proving that $\PSS$ is an isomorphism, an explicit inverse is constructed in the following. The adequate moduli space
$\M^{\mathrm{PSS, inv}}(\bar{x},q)$ is composed of pairs $(u,\gamma)$ such that
\beq
u:\R\times S^1\pf M\;\text{ with}\;E(u)<\infty\quad\text{and}\quad \gamma:[0,\infty)\pf M
\eeq
solve the equations
\beq
\partial_su+J(u)\big(\partial_tu-\beta(-s)X_{H}(t,u)\big)=0
        \quad\text{and}\quad\dot{\gamma}+\nabla^{g} f\circ\gamma=0\,,
\eeq
and fulfill boundary conditions $u(-\infty)=x$, $u(\infty)=\gamma(0)$ and $\gamma(+\infty)=q$,
where $q\in \Crit(f)$ is a critical point of $f$ and $\bar{x}\in\Pt(H)$. Finally, we impose the homotopy condition $\om(d_x\#u)=0=c_1(d_x\#u)$.
As before, for generic choices we obtain a smooth manifold $\M^{\mathrm{PSS, inv}}(\bar{x},q)$ of
$\dim \M^{\mathrm{PSS, inv}}(\bar{x},q)=m-\Morse(q)-\CZ(\bar{x})$, which is compact in dimension 0. We define
\bea
\PSS^{-1}:\CM^*(f)&\pf\CF_{m-*}(H)\\
q\;\;&\mapsto\!\!\!\!\sum_{\CZ(\bar{x})=m-\Morse(q)}\#_2\M^{\mathrm{PSS, inv}}(\bar{x},q)\cdot \bar{x}\,.
\eea
As above a suitable compactification of the 1-dimensional components of  $\M^{\mathrm{PSS, inv}}(\bar{x},q)$ entails
\beq
\PSS^{-1}:\HM^*(f)\pf\HF_{m-*}(H)\,.
\eeq
Up to here the nomenclature $\PSS^{-1}$ is purely formal. In what follows a justification is sketched.
The following arguments are taken from
\cite{Piunikhin_Salamon_Schwarz_Symplectic_Floer_Donaldson_theory_and_quantum_cohomology}.

Composing $\PSS^{-1}\circ\PSS:\HF_*(H)\pf\HF_*(H)$ amounts to
\beq
\PSS^{-1}\circ\PSS(\bar{x})=\sum_{\bar{y}}\sum_q\#_2\M^{\mathrm{PSS, inv}}(\bar{y},q)\cdot\#_2\M^{\mathrm{PSS}}(q,\bar{x})\cdot \bar{y}\,,
\eeq\\[-2ex]
i.e.~counting Floer cylinders emanating from $\bar{x}\in\Pt(H)$ connecting to gradient flow half-trajectories ending at some
critical point $q\in \Crit(f)$. From $q$ further gradient flow half-trajectories finally connect to
Floer cylinders ending at $\bar{y}\in\Pt(H)$. The idea is to form a cobordism between these configurations and the
identity as follows:
\begin{enumerate}
\item We glue the two gradient flow half-trajectories at $q$ to obtain a \textbf{finite length} gradient flow trajectory.
\item Shrink that finite length to zero, i.e.~we end up with two Floer cylinders, which
    meet at the same point.
\item Due to the cut-off functions on the respective ends, these two Floer cylinders are holomorphic near that point.
    We employ a glueing theorem for holomorphic curves to obtain one Floer cylinder passing from $\bar{x}$ to $\bar{y}$.
\item The Hamiltonian function for this Floer cylinder is not yet $H$ due to the cut-off. Therefore, we choose a
    (compactly supported) homotopy changing this Hamiltonian function to $H$.
\item We end up with counting Floer cylinders for the Hamiltonian function $H$ connecting periodic solutions $\bar{x}$ and $\bar{y}$
    of the \textbf{same} Conley-Zehnder index. For dimension reasons generically there are no such solutions unless $\bar{x}=\bar{y}$,
    in which case there is only the constant solution. This defines the identity homomorphism on
    Floer homology.
\end{enumerate}
For the other composition $\PSS\circ\PSS^{-1}:\HM_*(f)\pf\HM_*(f)$ we obtain
\beq
\PSS\circ\PSS^{-1}(q)=\sum_p\sum_{\bar{x}}\#_2\M^{\mathrm{PSS}}(p,\bar{x})\cdot\#_2\M^{\mathrm{PSS, inv}}(\bar{x},q)\cdot p\,.
\eeq\\[-2ex]
Here two Floer cylinders meet at a periodic solution $\bar{x}$. This leads to the following cobordism:
\begin{enumerate}
\item We glue the two Floer cylinders at $\bar{x}$ to obtain a sphere obeying Floer's equation.
\item Since the Hamiltonian term vanishes outside a neighborhood of the equator of this sphere we choose a homotopy
    of the Hamiltonian term to zero and end up with a holomorphic sphere.
\item Assuming for the moment that $\om|_{\pi_2(M)}=0$, this sphere is constant. We reduced the problem to count
    gradient trajectories for index difference 0. As above this defines the identity homomorphism.
\end{enumerate}
In the paper \cite{Piunikhin_Salamon_Schwarz_Symplectic_Floer_Donaldson_theory_and_quantum_cohomology} and in the book
\cite[chapter 12.1]{McDuff_Salamon_J_holomorphic_curves_and_symplectic_topology}
are series of figures picturise the above ideas.
We conclude that our notation is meaningful, i.e.~$\PSS$ is an isomorphism with inverse
\beq
(\PSS)^{-1}=\PSS^{-1}.
\eeq
in the case $\om|_{\pi_2(M)}=0$.
In the semi-positive case \eqref{property:semi_positive} step (3) becomes much more delicate since, in general, there will exist holomorphic
spheres. Then an elaborate transversality argument is necessary to allow for a time independent almost complex structure along the sphere. Since
only two points on the sphere are fixed by the gradient trajectories, an $S^1$-symmetry remains. In particular, solutions
come in 1-dimensional families as long as the sphere is non-constant. This contradicts the fact that the moduli
space that we started with has dimension 0. We end up with the same conclusion as in (3).
\begin{Def}\label{def:homologically_essential}
An element $\bar{x}\in\CF_*(H)$ is called \textit{essential} if there are $p,q\in\Crit(f)$ s.t.
\beq\label{property:homologically_essential}
\M^{\mathrm{PSS}}(p,\bar{x})\not=\emptyset\quad\text{and}\quad\M^{\mathrm{PSS, inv}}(\bar{x},q)\not=\emptyset\,.
\eeq
\end{Def}

\begin{Lemma}\label{Lemma_action_estimate_for_homologically_essential}
For essential elements $\bar{x}\in\CF_*(H)$ we can estimate (cf.~inequality \eqref{eqn:energy_of_plumbers_helper_solutions})
\beq
-\int_0^1\sup_MH\,dt\;\leq\;\A_{H}(\bar{x})\;\leq\;-\int_0^1\inf_MH\,dt\,.
\eeq
\end{Lemma}
\begin{Prop}\label{prop:non_vanishing_in_hommology_implies_essential}
All periodic orbits representing non-zero homology classes are essential.
\end{Prop}
\begin{proof}
Pick $\bar{x}\in\CF_*(H)$ such that $\partial{\bar{x}}=0$ and $[\bar{x}]\not=0$.
If $\M^{\mathrm{PSS}}(p,\bar{x})=\emptyset$ for all $p\in\Crit(f)$, then $\PSS([\bar{x}])=0$.

The existence of $q$ such that $\M^{\mathrm{PSS, inv}}(\bar{x},q)\not=\emptyset$ follows via \Poincare
duality, cf.\cite{Schwarz_Matthias_Action_spectrum}.
\end{proof}
\noindent\textbf{The action filtration and the definition of spectral capacities.}\quad
The action filtration and the action selector are used only to state and prove theorem \ref{thm:high_minimal_Maslov_number_implies_intersection}.
Therefore, we are very brief about the definitions and properties. Detailed expositions can be found in
\cite{Schwarz_Matthias_Action_spectrum,
McDuff_Salamon_J_holomorphic_curves_and_symplectic_topology}.

We define the \textit{action filtration} on Floer homology by
\beq
\CF_k^a(H):=\big\{\!\!\sum_{\CZ(\bar{x})=k}a_{\bar{x}}\bar{x}\mid a_{\bar{x}}=0\text{ for }\A_H(\bar{x})>a\big\},
\eeq
and check that $\big(\CF_k^a(H),\partial\big)$ forms a subcomplex of the Floer complex
(cf.~\cite{McDuff_Salamon_J_holomorphic_curves_and_symplectic_topology}) due to the chosen Novikov condition. In particular,
we obtain the long exact sequence
\beq
\cdots\pf\HF_*^a(H)\stackrel{i_H^a}{\pf}\HF_*(H)\stackrel{j_H^a}{\pf}\HF_*^{(a,\infty]}(H)\pf\cdots\;.
\eeq
This gives rise to the \textit{action selector}
\beq
c(\alpha,H):=\inf\{a\in\R\mid \PSS^{-1}(\alpha)\in\im i_H^a\}
\eeq
for all $0\not=\alpha\in\QH^*(M)$. We note here, that the \textit{quantum cohomology} of the symplectic manifold comes into play and
refer the reader to \cite[chapter 11.1]{McDuff_Salamon_J_holomorphic_curves_and_symplectic_topology}.

The action selector is well-defined, i.e.~$c(\alpha,H)$ is a finite number,
since each representation of the class $\PSS^{-1}(\alpha)$ contains
an element $\bar{x}\in\Pt(H)$ with maximal action value, due to the Novikov condition (but there might not
be such an element with minimal action value).
\\[1ex]\textbf{Properties of the action selector.}\quad
A full list of properties of the action selector $c(\alpha,H)$ can be found in
\cite[chapter 12.4]{McDuff_Salamon_J_holomorphic_curves_and_symplectic_topology}. We shall need the following:
For $\alpha, \beta\in\QH^*(M)$, $A\in\Gamma$ and Hamiltonian functions $H,K:S^1\times M\pf\R$ we have\\[1ex]
\begin{tabular}{ll}
(\textit{Zero}) \hspace{15ex} & $c(\alpha,0)=0$\\[1ex]
(\textit{Novikov action}) & $c(\alpha\ast e^A,H)=c(\alpha,H)-\om(A)$\\[1ex]
(\textit{Product}) & $c(\alpha\ast\beta,H\#K)\leq c(\alpha,H)+c(\beta,K)$\\[1ex]
(\textit{\Poincare duality}) & $c(1,H)\geq -c([\om^m],H^{(-1)})$\\[1ex]
(\textit{Non-degeneracy}) & for each $H\not\equiv0$ there exists a $\delta>0$ such that\\[0.5ex]
& $c(1,H)+c(1,H^{(-1)})\;\geq\; c(1,H)-c([\om^m],H)\;\geq\;\delta\,.$
\end{tabular}\\[1ex]
Here, $H\#K(t,m)=H(t,m)+K\big(t,(\psi_H^t)^{-1}(m)\big)$ and $\frac{d}{dt}\psi_H^t=X_H(t,\psi_H^t)$;
furthermore, we set $H^{(-1)}(t,m)=-H(-t,m)$ and $\ast$ denotes the \textit{quantum-cup-product} in $\QH^*(M)$.

With the help of the action selector we define two norms (cf.~\cite{McDuff_Salamon_J_holomorphic_curves_and_symplectic_topology,
Schwarz_Matthias_Action_spectrum})
\beq
\gamma(H):=c(1,H)-c([\om^m],H)\zw\text{and}\zw\tilde{\gamma}(H):=c(1,H)+c(1,H^{(-1)})\,.
\eeq
\textit{Non-degeneracy} implies $\tilde{\gamma}(H)\geq\gamma(H)$. These norms give rise to two symplectic capacities as follows. For
an open subset $U\subset M$ we set
\beqn
c_\gamma(U)=\sup\{\gamma(H)\mid \mathrm{supp}X_H\subset S^1\times U\}\zw\text{and}\zw
\tilde{c}_{\tilde{\gamma}}(U)=\sup\{\tilde{\gamma}(H)\mid \mathrm{supp}X_H\subset S^1\times U\}\,.\\
\eeq
Again \textit{non-degeneracy} implies $\tilde{c}_{\tilde{\gamma}}(U)\geq c_\gamma(U)$.
We call the capacity $ c_\gamma(U)$ the \textit{spectral capacities} of $U$. Most interesting to us is the case $U=M$. For a Lagrangian submanifold
$L\subset M$ we define $c_\gamma(L)$ as $c_\gamma(L):=\inf\{c_\gamma(U)\mid U\text{ open},\;L\subset U\}$
and $\tilde{c}_{\tilde{\gamma}}(L)$ accordingly.
\section{Basic theorems}
\label{sec:basic_theorems}

The basic theorems which all other results rely on are the following theorems \ref{thm:basic_theorem} and \ref{thm:basic_theorem_with_low_energy}.

\begin{Thm}\label{thm:basic_theorem}
Let $(M^{2m},\om)$ be a monotone closed symplectic manifold and $L\subset M$ a
monotone closed Lagrangian submanifold with $N_L\geq 2$ (cf.~definition \ref{def:minimal_numbers}).
For each closed submanifold $S\subset L$ such that $s:=\dim S$ satisfies
\beq
s>\dim L+1-N_L\,,\qquad\text{i.e.}\zw N_L>\mathrm{codim}_L S+1\,,
\eeq
the image of $[S]\in\H_s(M;\Z/2)$ under the $\PSS$-isomorphism is represented by the cycle
\beq\label{eqn:Floer_representation_PSS^(-1)([S])}
\PSS^{-1}\big([S]\big)=\left[\sum_{\bar{x}\in\Pt_{s-2m}(H)}
        \#_2\left(\M^+_{(L,S)}(\bar{x};J,H;\mathbf{0})\right)\cdot \bar{x}\right]\in\HF^{2m-s}(H)\,,
\eeq
in the Floer cohomology of the (generically chosen) Hamiltonian $H:S^1\times M\rightarrow\R$, and dually
\beq\label{eqn:Floer_representation_PSS^(-1)(PD[S])}
\PSS^{-1}\big(\PD[S]\big)=\left[\sum_{\bar{x}\in\Pt_{s-2m}(H)}
        \#_2\left(\M^-_{(L,S)}(\bar{x};J,H;\mathbf{0})\right)\cdot \bar{x}\right]\in\HF_{s-2m}(H)\,.
\eeq
The moduli spaces $\M^\pm_{(L,S)}(\bar{x};J,H;\mathbf{0})$ are defined with the help of a generically chosen
almost complex structures $J$ as
\beq\label{def:M_(L,S)^+}
\M^+_{(L,S)}(\bar{x};J,H;\mathbf{0}):=\left\{u:[0,\infty)\times S^1\pf(M,L)\,\left|\,
    \begin{aligned}
        &\db u=0,\;&&[u\#\bar{d}_x]=\mathbf{0}\\[0.5ex]
        &u(+\infty)=x,&&u(0,0)\in S
    \end{aligned}
\,\right.\right\}\,
\eeq
respectively
\beq\label{def:M_(L,S)^-}
\M^-_{(L,S)}(\bar{x};J,H;\mathbf{0}):=\left\{u:(-\infty,0]\times S^1\pf(M,L)\,\left|\,
    \begin{aligned}
        &\db u=0,\;&&[d_x\#u]=\mathbf{0}\\[0.5ex]
        &u(-\infty)=x,&&u(0,0)\in S
    \end{aligned}
\,\right.\right\}\,,
\eeq
where $\bar{\partial}_{J,H}$ abbreviates Floer's equation, namely $\partial_su+J(u)\big(\partial_tu-X_H(t,u)\big)=0$,
and $u:[0,\infty)\times S^1\pf(M,L)$ means $u(\{0\}\times S^1)\subset L$.
Recall that $\bar{x}=[x,d_x]$. The homotopical condition $[d_x\#u]=\mathbf{0}$ in the definition of
$\M^-_{(L,S)}(\bar{x};J,H;\mathbf{0})$ is meant in the Novikov sense, i.e.~$\om([d_x\#u])=0=\Mas([d_x\#u])$
and accordingly for $\M^+_{(L,S)}(\bar{x};J,H;\mathbf{0})$.
\end{Thm}

We count Floer \textit{half-cylinders} $u$ which are asymptotic to a periodic orbit and have boundary on the Lagrangian submanifold,
i.e.~$u(0,t)\in L$ for all $t\in S^1$.
An important feature of the moduli spaces is the \emph{marking} of the half-cylinders: we
require them to map the marked point $(0,0)$ on the half-cylinder to $S$.
The closed submanifold $S$ can be replaced by any singular chain representing a cycle on $L$.

\begin{Rmk}\label{rmk:extremal_cases_of_basic_theorems}
If $S=L$, we have $\dim L>\dim L+1-N_L$, since $N_L\geq2$. In particular,
theorem \ref{thm:basic_theorem} implies that we always can represent the class $[L]\in\H_m(M;\Z/2)$.

If $S=\mathrm{pt}$ and $N_L>\dim L+1$, we obtain a Floer-theoretic representation of the class $[\mathrm{pt}]\in\H_0(M;\Z/2)$ and its \Poincare dual
$[\om^m]\in\H^{2m}(M;\Z/2)$.
The condition $N_L>\dim L+1$ excludes all \emph{displaceable} (cf.~definition \ref{def:displaceable}) Lagrangian submanifolds. This is no coincidence in
view of theorem \ref{thm:displaceable_implies_empty_moduli_spaces}.
\end{Rmk}
\begin{Thm}\label{thm:basic_theorem_with_low_energy}
If we assume in addition to the assumptions of theorem \ref{thm:basic_theorem}
that the Hamiltonian function $H:S^1\times M\pf\R$ satisfies
\beq
||H||\leq\min\big\{A_M,\,A_L\big\}
\eeq
(see definition \ref{def:minimal_numbers} for $A_M,A_L$), then the representations \eqref{eqn:Floer_representation_PSS^(-1)([S])}
and \eqref{eqn:Floer_representation_PSS^(-1)(PD[S])} hold for all closed submanifolds $S\subset L$ regardless of their dimension.
\end{Thm}

The basic theorems \ref{thm:basic_theorem} and \ref{thm:basic_theorem_with_low_energy} deal with the \textit{extrinsic
topology} of the Lagrangian submanifold: They provide a Floer theoretic representation of the homomorphism
$\iota:\H_*(L;\Z/2)\pf\H_*(M;\Z/2)$ induced by the inclusion $\iota:L\subset M$.
\section{Proof of the basic theorems}
\label{sec:Proof_of_the_basic_theorems}

The proofs of theorems \ref{thm:basic_theorem} and \ref{thm:basic_theorem_with_low_energy} rely
both on the compactness theorem \ref{thm:compactness} stated below
for the moduli spaces $\M^\pm_{(L,S)}(\bar{x};J,H;\mathbf{0})$ in dimensions 0 and 1. The respective assumption
on the dimension of $S$ or the Hofer norm of $H$ prohibits bubbling-off,
which in turn results in the compactness of the moduli spaces in question.
Before proving this compactness result we make sure that the moduli spaces are smooth manifolds for generic
choice of $H$ and $J$.
\begin{Thm}\label{thm:transversality}
Let $L$ be a closed Lagrangian submanifold of the closed symplectic manifold $(M,\om)$ and $S\subset L$ a closed
submanifold.
For a generic Hamiltonian function $H:S^1\times M\rightarrow\R$ and a generic
almost complex structure $J$ the moduli spaces
\beq
\M^-_{(L,S)}(\bar{x};J,H;\mathbf{0})\quad\text{and}\quad
\M^+_{(L,S)}(\bar{x};J,H;\mathbf{0})
\eeq
defined in \eqref{def:M_(L,S)^+} and \eqref{def:M_(L,S)^-} are smooth manifolds of dimension
\beq
\dim\M^\pm_{(L,S)}(\bar{x};J,H;\mathbf{0})=\pm\CZ(\bar{x})-\dim L+\dim S\,,
\eeq
where $\CZ(\bar{x})$ denotes the Conley-Zehnder index of the periodic orbit $\bar{x}\in\Pt(H)$.
\end{Thm}

\begin{proof}
This is a standard result in Floer theory (cf.~\cite{Floer_Hofer_Salamon_Transversality_in_elliptic_Morse_theory_for_the_symplectic_action})
with one minor modification, cf.~remark \ref{rmk:on_transversality}\,(i).
First of all, for a non-degenerate $H$ we can regard $\db=\partial_s+J(\partial_t-X_H)$ as a Fredholm-section in a
suitable Banach-bundle and identify the moduli spaces with the vanishing locus of this Fredholm-section. For generic $H$ and $J$,
this section will be transverse to the zero-section. In particular, by the implicit function theorem,
the moduli spaces are smooth finite-dimensional manifolds. Computing the dimension of the moduli
spaces is achieved by the Riemann-Roch theorem and additivity of the Fredholm index. All details within the setting of this paper can
be found in \cite{Albers_PhD}.
\end{proof}

\begin{Rmk}\label{rmk:on_transversality}
\begin{enumerate}\renewcommand{\labelenumi}{(\roman{enumi})}
\item To prove the transversality result we need to assume, in addition to $H$ being non-degenerate, that there are no periodic
orbits of $H$ lying entirely on $L$. This excludes $s$-independent solutions of the Floer
equation, for which we were not able to prove the necessary transversality statement. In contrast to the construction
of Floer homology there is \textbf{no} (immediate) \textbf{automatic transversality} result for half-cylinders. Anyway,
this additional assumption on $H$ is clearly fulfilled generically.
\item The moduli space $\M^\pm_{(L,S)}(\bar{x};J,H;\mathbf{0})$ is a submanifold of
$\M^\pm_{(L,L)}(\bar{x};J,H;\mathbf{0})$ for generic almost complex structures, because the evaluation map
\beq
\ev:\M^\pm_{(L,L)}(\bar{x};J,H;\mathbf{0})\pf L\qquad\ev(u):=u(0,0)
\eeq
is a submersion if regarded on a universal moduli space, cf.~\cite[chapter 3.4]{McDuff_Salamon_J_holomorphic_curves_and_symplectic_topology},
and
\beq
\M^\pm_{(L,S)}(\bar{x};J,H;\mathbf{0})=\ev^{-1}(S)\,.
\eeq
\end{enumerate}
\end{Rmk}\noindent
To complete the proofs of theorems \ref{thm:basic_theorem} and \ref{thm:basic_theorem_with_low_energy} we
need to show that

\begin{enumerate}\renewcommand{\labelenumi}{\alph{enumi})}
\item the 0-dimensional moduli spaces are compact,
\item the chains $\sum_{\bar{x}}
        \#_2\M^\pm_{(L,S)}(\bar{x};J,H;\mathbf{0})\cdot \bar{x}$ are cycles in the Floer complex,
\item they represent the claimed (co)-homology classes.
\end{enumerate}

The first two points are the content of the following compactness results, theorem \ref{thm:compactness}.
Since the assertions of theorem \ref{thm:basic_theorem_with_low_energy} are homological it suffices in the second case of theorem \ref{thm:compactness}
to consider only essential elements since all elements representing non-vanishing homology classes are essential,
cf.~proposition \ref{prop:non_vanishing_in_hommology_implies_essential}.

\begin{figure}[ht]
\begin{center}
\input{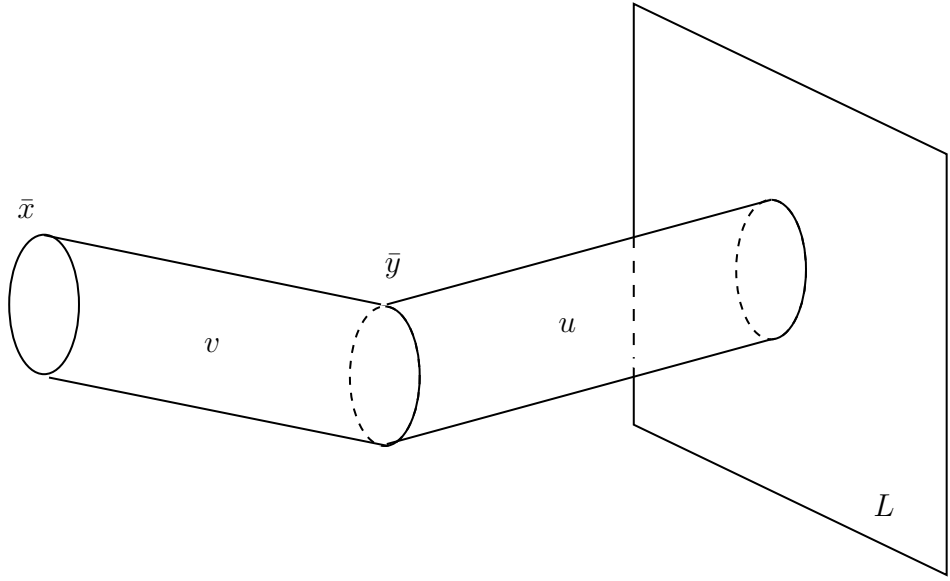}
\caption{A broken solution representing $(v,u)\in\partial\M^-_{(L,S)}(\bar{x};J,H;\mathbf{0})$.}\label{fig_half_cylinder_breaking}
\end{center}
\end{figure}

\begin{Thm}\label{thm:compactness}
Let $(M^{2m},\om)$ be a closed monotone symplectic manifold and $L\subset M$ a
monotone closed Lagrangian submanifold with $N_L\geq2$.
If we choose a closed submanifold $S\subset L$ and a Hamiltonian function $H:S^1\times M\rightarrow\R$
such that one of the following holds:
\begin{enumerate}\renewcommand{\labelenumi}{$\bullet$}
\item \quad$\dim S>\dim L+1-N_L$\,,\\[-0.5ex]
\item \quad$||H||\leq\min\big\{A_M,\,A_L\big\}$ and $\bar{x}\in\Pt(H)$ is essential (cf.~definition \ref{def:homologically_essential}),
\end{enumerate}
then the 0-dimensional moduli spaces $\M^\pm_{(L,S)}(\bar{x};J,H;\mathbf{0})$ are compact.

Furthermore, the 1-dimensional moduli spaces are compact up to splitting-off one Floer cylinder,
i.e.~they can be compactified in such a way
that the boundary of the compactification, which we denote by the same symbol, decomposes as follows
(cf.~formula \eqref{eqn:moduli_space_Floer_cylinders})
\bea\label{decomposition_brdy_M_L(y)=Mh(y,y')xM_l(y')}
\partial\M^-_{(L,S)}(\bar{x};J,H;\mathbf{0})
        &=\bigcup_{\CZ(\bar{y})=\CZ(\bar{x})+1}\Mh(\bar{x},\bar{y};J,H)\times\M^-_{(L,S)}(\bar{y};J,H;\mathbf{0}),\\[2ex]
\partial\M^+_{(L,S)}(\bar{x};J,H;\mathbf{0})
        &=\bigcup_{\CZ(\bar{y})=\CZ(\bar{x})-1}\M^+_{(L,S)}(\bar{y};J,H;\mathbf{0})\times\Mh(\bar{y},\bar{x};J,H).
\eea
A typical element in $\partial\M^-_{(L,S)}(\bar{x};J,H;\mathbf{0})$ is depicted in figure \ref{fig_half_cylinder_breaking}.
\end{Thm}
\begin{proof}
The following inequality for elements $u\in\M^-_{(L,S)}(\bar{x};J,H;\mathbf{0})$\\[-2ex]
\beq\label{eqn:energy_inequality}
E(u)=\int_{-\infty}^0\int_0^1|\partial_s u|^2dt\,ds\leq \om([d_x\#u])-\A_H(\bar{x})-\int_0^1\inf_LH(t,\cdot)dt\,.
\eeq
is easily derived using Floer's equation. Let us abbreviate from now on
\beq
\inf_LH:=\int_0^1\inf_LH(t,\cdot)dt\qquad\text{and}\qquad\sup_LH:=\int_0^1\sup_LH(t,\cdot)dt\,.
\eeq
By the definition of $\M^-_{(L,S)}(\bar{x};J,H;\mathbf{0})$ we have $\om([d_x\#u])=0$, therefore \eqref{eqn:energy_inequality} reads
\beq\label{eqn:energy_estimate_M_(L,S)^-}
0\leq E(u)\leq-\A_H(\bar{x})-\inf_LH
\eeq
and analogously for elements $u\in\M^+_{(L,S)}(\bar{x};J,H;\mathbf{0})$
\beq\label{eqn:energy_estimate_M_(L,S)^+}
0\leq E(u)\leq+\A_H(\bar{x})+\sup_LH\,.
\eeq
Hence, the energy of elements of the moduli spaces is uniformly bounded.\\[1.5ex]
\textbf{The case $S=L$.}$ $\\[0.5ex]
Since the energy is uniformly bounded and the marking $u(0,0)\in S=L$ for an element
in the moduli space is obsolete we conclude by the standard compactness results in Floer theory that
sequences $(u_n)\subset\M^-_{(L,L)}(\bar{x};J,H;\mathbf{0})$ converge in the Gromov-Hausdorff topology
\beq
u_n\rightarrow (u_\infty\,;\,v_1,\ldots,v_\Gamma\,;\,s_1,\ldots,s_\Sigma\,;\,d_1,\ldots,d_\Delta)
\eeq
where
\begin{itemize}
\item $u_\infty\in\M^\pm_{(L,L)}(\bar{x}_0;J,H;\mathbf{0})$,\\[-1ex]
\item $v_\gamma\in\Mh(\bar{x}_{\gamma},\bar{x}_{\gamma-1};J,H)$, where $\bar{x}_\gamma\in\Pt(H)$
        and $\bar{x}_{\Gamma}=\bar{x}$,\\[-1ex]
\item $\{s_\sigma\}$ are holomorphic spheres,\\[-1ex]
\item $\{d_\delta\}$ are holomorphic disks,\\[-1ex]
\end{itemize}
i.e.~sequences converge to a finite family of adjacent Floer cylinders starting at $\bar{x}$ and finally connecting
to a half-cylinder in $\M^-_{(L,L)}(\bar{x}_0;J,H;\mathbf{0})$. Furthermore, there are finitely many
holomorphic spheres and disks attached to this configuration.

The Fredholm index behaves additively \wrt this convergence, i.e.~the Fredholm index at an element $u_n$,
which equals the dimension of $\M^-_{(L,L)}(\bar{x};J,H;\mathbf{0})$, is
\beq\label{eqn:Fredholm_index_formula_bubbling}
\ind\F_{u_n}=\ind\F_{u_\infty}+\sum_\gamma\ind\F_{v_\gamma}+\sum_\sigma2c_1(s_\sigma)+\sum_\delta\Mas(d_\delta)\,.
\eeq
Since the Lagrangian submanifold and therefore the symplectic manifold is monotone, the Chern class and the Maslov index
evaluate positively in case that $s_\sigma$ and $d_\delta$ are non-constant.
Moreover, by assumption the Maslov index is at least 2,
and if the Floer cylinders $v_\gamma$ are $s$-dependent, the Fredholm indices are at least 1.
In particular, if the moduli space $\M^-_{(L,L)}(\bar{x};J,H;\mathbf{0})$
has dimension 0 all maps but $u_\infty$ must be trivial. In other words, the moduli space
$\M^-_{(L,L)}(\bar{x};J,H;\mathbf{0})$ is compact.

If the moduli space $\M^-_{(L,L)}(\bar{x};J,H;\mathbf{0})$ is 1-dimensional, then a sequence converges to the limit half-cylinder $u_\infty$
and at most one further Floer cylinder $v_\gamma$. In particular, at most one Floer cylinder may split-off. This gives rise to the compactification
with the decomposition of the boundary as asserted in the theorem.

Since twice the first Chern class evaluated on a sphere is a multiple of the minimal Maslov number $N_L$
(cf.~remark \ref{rmk:orientable_Lagrangian_minimal_Maslov_geq_2}) we obtain the following refined statement.
If the dimension of the moduli space is less than the minimal Maslov number
\beq
\dim \M^-_{(L,L)}(\bar{x};J,H;\mathbf{0})=-\CZ(\bar{x})<N_L
\eeq
the moduli space $\M^-_{(L,L)}(\bar{x};J,H;\mathbf{0})$ is compact up to splitting-off a finite number of Floer cylinders.
No holomorphic spheres or disks bubble off. A completely analogous statement holds for the moduli space $\M^+_{(L,L)}(\bar{x};J,H;\mathbf{0})$. \\[1.5ex]
\textbf{The case $\dim S>\dim L+1-N_L$.}$ $\\[0.5ex]
We already noted above that the marked moduli space $\M^\pm_{(L,S)}(\bar{x};J,H;\mathbf{0})$ is a submanifold of the
moduli space $\M^\pm_{(L,L)}(\bar{x};J,H;\mathbf{0})$ for generic
almost complex structures.

If $\M^\pm_{(L,S)}(\bar{x};J,H;\mathbf{0})$ is 0 or 1-dimensional and in addition we require the condition $\dim S>\dim L+1-N_L$ to hold,
we conclude that $\dim \M^\pm_{(L,L)}(\bar{x};J,H;\mathbf{0})<N_L$.
As explained at the end of the first case we obtain the statement of theorem \ref{thm:compactness} in this case.\\[1.5ex]
\textbf{The case $||H||\leq\min\big\{A_M,\,A_L\big\}$ and $\bar{x}\in\Pt(H)$ is essential.}$ $\\[0.5ex]
Without the dimension restriction on the submanifold $S$ we cannot expect compactness of all 0 and 1-dimensional moduli spaces
$\M^\pm_{(L,S)}(\bar{x};J,H;\mathbf{0})$. In fact, some of our applications are derived with help of this non-compactness.
But as long as all elements in $\M^\pm_{(L,S)}(\bar{x};J,H;\mathbf{0})$ have energy
less than the minimal energy of a non-constant holomorphic sphere or disk, a bubble would take away more energy than is available.

We assume that the Hamiltonian function $H$ has sufficiently small Hofer norm, namely $||H||\leq\min\big\{A_M,\,A_L\big\}$,
and the periodic orbit $\bar{x}$ is essential, cf.~definition \ref{def:homologically_essential}.
For such a periodic orbit $\bar{x}$ we recall from lemma
\ref{Lemma_action_estimate_for_homologically_essential} the inequalities $-\sup_MH\leq\A_H(\bar{x})\leq-\inf_MH$.
Combining this with inequality \eqref{eqn:energy_estimate_M_(L,S)^-} resp.~\eqref{eqn:energy_estimate_M_(L,S)^+}
we obtain
\beq\label{eqn:inequality_in_proof_of_compacness_theorem}
E(u)\leq-\A_H(\bar{x})-\inf_LH\leq\sup_MH-\inf_MH=||H||
\eeq
for an element $u\in\M^-_{(L,S)}(\bar{x};J,H;\mathbf{0})$ and analogously
\beq
E(u)\leq+\A_H(\bar{x})+\sup_LH\leq-\inf_ML+\sup_MH=||H||
\eeq
for an element $u\in\M^+_{(L,S)}(\bar{x};J,H;\mathbf{0})$.

In particular, bubbling-off cannot occur for a sequence $(u_n)\subset\M^\pm_{(L,S)}(\bar{x};J,H;\mathbf{0})$.
Indeed, the energy $E(u_\infty)$ of the limit solution $u_\infty$ is
positive by the assumption on $H$ that no periodic orbits lie entirely on $L$. Therefore,
\beq\label{eqn:final_inequality_in_proof_of_compacness_theorem}
E(b)<E(u_\infty)+E(b)\leq \liminf E(u_n)\leq||H||\leq\min\big\{A_M,\,A_L\big\}\,.
\eeq
We obtain $E(b)=0$. Since we excluded bubbling-off of holomorphic spheres and disks the relevant moduli spaces are compact or can be
compactfied as asserted. This finishes the proof of theorem \ref{thm:compactness}.
\end{proof}

\begin{Rmk}\label{rmk:compactness_with_direct_inequality}
For later purposes we remark that in the proof of theorem \ref{thm:compactness} the assumption $||H||\leq\min\{A_M,A_L\}$ is used only in
combination with inequality \eqref{eqn:inequality_in_proof_of_compacness_theorem} to derive the crucial inequality $E(u)\leq\min\{A_M,A_L\}$.
Therefore, if we assume this inequality right away, the assertion of theorem \ref{thm:compactness} still
holds.
\end{Rmk}

\begin{Rmk}
A geometric explanation for the hypothesis $\dim S>\dim L+1-N_L$ is that a holomorphic disk bubbling-off may
\textit{take away the marking}.
If we disregard the marking in a configuration as depicted in figure \ref{figure_half_cylinder_with_marking_and_bubbling}, it lies in the boundary
of the higher dimensional moduli space $\M^\pm_{(L,L)}(\bar{x};J,H;\mathbf{0})$. The assumption $\dim S>\dim L+1-N_L$ excludes
such configurations in $\M^\pm_{(L,L)}(\bar{x};J,H;\mathbf{0})$ for index reasons so that no bubble can take away the marking.
\begin{figure}[ht]
\begin{center}
\input{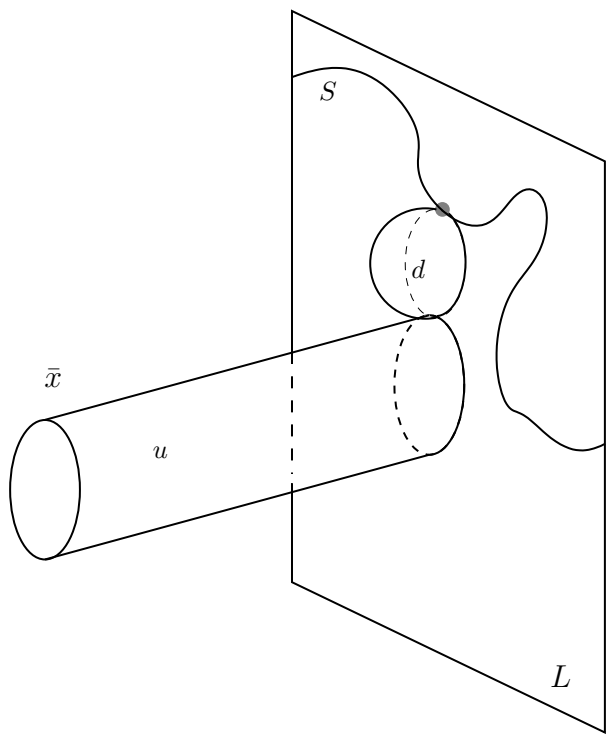}
\caption{A configuration in which the bubble $d$ ''takes away'' the marked point.}
\label{figure_half_cylinder_with_marking_and_bubbling}
\end{center}
\end{figure}
\end{Rmk}
Theorem \ref{thm:compactness} immediately implies that the chain defined in
theorem \ref{thm:basic_theorem} resp.~\ref{thm:basic_theorem_with_low_energy} is well-defined and in fact a cycle.
To finish the proof of theorems \ref{thm:basic_theorem} and \ref{thm:basic_theorem_with_low_energy} we need
to show that the $\PSS$-isomorphism maps $[S]\in\H_s(M;\Z/2)$ of the homology class of this cycle.

\begin{proof}[\textsc{End of the proof of theorems \ref{thm:basic_theorem} and \ref{thm:basic_theorem_with_low_energy}}]\quad
We will define a cobordism relating the counting procedures defining $\PSS^{-1}([S])$ and the explicitly given cycle.
This shows that they are chain homotopic and therefore agree in homology.

We define the moduli space $\M^{\PSS,-}_{(L,S)}(q;g,f;J,H;\mathbf{0})$ as the set of triples $(R,\gamma,u)$, where
\beq
R\geq0,\quad\gamma:(-\infty,0]\pf M,\quad u:(-\infty,R]\times S^1\pf M
\eeq
subject to the equations
\beq
\dot{\gamma}+\nabla^{g} f\circ\gamma=0
        \quad\text{and}\quad\partial_su+J(u)\big(\partial_tu-\beta(s)X_{H}(t,u)\big)=0\,,
\eeq
where $\beta:\R\rightarrow[0,1]$ is a smooth cut-off function satisfying $\beta(s)=0$ for $s\leq\tfrac12$ and $\beta(s)=1$ for $s\geq1$.
The function $f:M\pf\R$ is an auxiliary Morse function and $\nabla^{g}$
is the gradient \wrt an auxiliary Riemannian metric $g$ on M.
Furthermore, $(R,\gamma,u)\in\M^{\PSS,-}_{(L,S)}(q;g,f;J,H;\mathbf{0})$ has to satisfy the following boundary conditions
\beq
E(u)<\infty\,,\;\;\gamma(-\infty)=q\,,\;\;\gamma(0)=u(-\infty)\;\;\text{and}\;\; u(R,t)\in L\;\forall t\in S^1\,,\;\; u(R,0)\in S\,,
\eeq
where $q\in\Crit(f)$ is a critical point of $f$. As in the definition of the $\PSS$-isomorphism the
finiteness of the energy of $u$ in conjunction with the cut-off function $\beta$ allows
for a continuous extension $u(-\infty)$. In particular, topologically we can think of $u$ as a disk.
Finally, we impose the homotopy condition $[u]=0$ for the disk $u$, i.e.~$\om(u)=0=\Mas(u)$.

The geometric ideas behind the definition of $\M^{\PSS,-}_{(L,S)}(q;g,f;J,H;\mathbf{0})$ are:
\begin{enumerate}
\item The cycle $\PSS\Big(\sum_{\bar{x}\in\Pt(H)}\#_2\big(\M^-_{(L,S)}(\bar{x};J,H;\mathbf{0})\big)\cdot \bar{x}\Big)$ is
defined by counting gradient half trajectories attached to Floer cylinders, i.e.~plumber's helper solutions, ending
at a periodic orbit $\bar{x}$, from where Floer half-cylinders start which end on $L$.
\item We glue such half-cylinders and Floer cylinders and obtain new half-cylinders which obey Floer's equation \wrt a Hamiltonian function
different from zero only near $L$.
\item We regard this Hamiltonian function as a compactly supported deformation of the zero Hamiltonian. Therefore, after a homotopy, we deal with
holomorphic cylinders which in fact are holomorphic disks since they have finite energy.
\item The assumption on the homotopy type of this disk immediately implies that it is constant. Thus, we are left with
counting gradient half-trajectories starting at a critical point and ending at $S$, due to the marking. It is easy to see
that this represents the class $\PD[S]\in\HM^s(M)$ in the Morse cohomology of $M$.
\end{enumerate}
Let us come back to our object of study: $\M^{\PSS,-}_{(L,S)}(q;g,f;J,H;\mathbf{0})$.
Combining the methods from the definition of the $\PSS$-isomorphism and from this chapter, it follows that
these spaces are smooth manifolds for generic choices of $g$, $J$ and $H$ and of dimension
\bea
\dim \M^{\PSS,-}_{(L,S)}(q;g,f;J,H;\mathbf{0})&=\dim S-\dim L-\big(\tfrac12\dim M-\Morse(q;f)\big)+1\\
    &=\dim S-\dim M+\Morse(q;f)+1\,.
\eea
The $+1$ accounts for the variable $R$. A straight forward computation implies for an element $(R,\gamma,u)\in\M^{\PSS,-}_{(L,S)}(q;g,f;J,H;\mathbf{0})$
that\\[-1.5ex]
\beq\label{eqn:energy_estimate_M_(L,S)^PSS,-}
E(u)\leq\beta(R)\cdot\Big(\sup_MH-\inf_LH\Big)\leq||H||\,.
\eeq
We note, that for $R=0$ this inequality implies $E(u)=0$.

Exactly the same arguments as in theorem \ref{thm:compactness} imply that the moduli spaces are compact in dimension 0 and can be
compactified in dimension 1. We denote again the compactification by the same letter. The boundary of the compactification decomposes as
\bean
\partial\M^{\PSS,-}_{(L,S)}(q;g,f;J,H;\mathbf{0})=
    &\bigcup_{
        \substack{
        q'\in\Crit(f)\\
        \Morse(q')=\Morse(q)-1}
    }\!\!\!\!\!\!\!\!\!\Mh^{\,\mathrm{Morse}}(q,q';f)\times\M_{(L,S)}^{\PSS,-}(q';g,f;J,H;\mathbf{0}) \\[1.5ex]
    \cup&\!\!\!\quad\bigcup_{
        \substack{
        \bar{x}\in\Pt(H)\\
        \CZ(\bar{x})=m-\Morse(q)}
    }\!\!\!\!\!\!\!\M^{\mathrm{PSS}}(q,\bar{x})\times\M_{(L,S)}^-(\bar{x};J,H;\mathbf{0})\\[1.5ex]
    \cup&\quad\Big\{(R,\gamma,u)\in\M_{(L,S)}^{\PSS,-}(q;g,f;J,H;\mathbf{0})\;\big|\; R=0\Big\}\;.
\eea
If we consider a sequence $(R_n,\gamma_n,u_n)_{n\in\N}$, the first union collates breaking of the gradient half-trajectory. In this case the sequence
$(R_n)$ converges. If the sequence $R_n\rightarrow\infty$, the sequence of half-cylinders $u_n$ breaks into a pair consisting
of a half-cylinder and a plumber's helper solution due to the chosen cut-off function $\beta$. This makes up the second union. The last union appears
for obvious reasons.

Since the moduli spaces $\M^{\PSS,-}_{(L,S)}(q;g,f;J,H;\mathbf{0})$ are compact in dimension 0, we can define 
\beq
\theta(S):=\!\!\!\!\!\!\!\!\!\!\sum_{\Morse(q)=\dim M-\dim S-1}\!\!\!\!\!\!\!\!\!\!
                        \#_2\M^{\PSS,-}_{(L,S)}(q;g,f;J,H;\mathbf{0})\cdot q\;\in\CM^*(f;\Z/2)\,,
\eeq
which is \emph{not} a cycle but merely a chain. Furthermore, we abbreviate the set
\beq
\M^{\PSS,-}_{(L,S)}(q;R=0):=\Big\{(R,\gamma,u)\in\M_{(L,S)}^{\PSS,-}(q;g,f;J,H;\mathbf{0})\;\big|\;R=0\Big\}
\eeq
and define the cycles
\begin{gather}
\rho(S):=\!\!\!\!\!\!\!\!\!\!\sum_{\Morse(q)=\dim M-\dim S}\!\!\!\!\!\!\!\!\!\!\#_2\M^{\PSS,-}_{(L,S)}(q;R=0)\cdot q\;\in\CM^*(f;\Z/2)\\[1ex]
\text{and}\qquad\Phi(S):=\sum_{\bar{x}\in\Pt_{s-2m}(H)}\#_2\left(\M^-_{(L,S)}(\bar{x};J,H;\mathbf{0})\right)\cdot \bar{x}\;\in\CF_*(H)\,.
\end{gather}
From the decomposition of moduli space $\partial\M^{\PSS,-}_{(L,S)}(q;g,f;J,H;\mathbf{0})$ we conclude
\beq
\delta^{\mathrm{Morse}}\big(\theta(S)\big)=\PSS\big(\Phi(S)\big)-\rho(S)\,,
\eeq
where the minus sign is arbitrary as long as we are working with $\Z/2$-coefficients.
We conclude that in cohomology $[\Phi(S)]=\PSS^{-1}\big([\rho(S)]\big)$. To finish the proof of theorems \ref{thm:basic_theorem}
and \ref{thm:basic_theorem_with_low_energy} we need to show $\rho(S)=\PD[S]$.

Either from equation \eqref{eqn:energy_estimate_M_(L,S)^PSS,-} or directly from the definition we conclude that for an element
$(0,\gamma,u)\in\M^{\PSS,-}_{(L,S)}(q;R=0)$, the map $u$ is constant.
Namely, since $R=0$ the Hamiltonian term in the Floer equation is 0, i.e.~$u$ is a holomorphic disk satisfying $E(u)=0$ due to the assumption $[u]=0$.
We conclude that $\M^{\PSS,-}_{(L,S)}(q;R=0)$ is in bijection to the space containing gradient half-trajectories $\gamma:(-\infty,0]\pf M$
such that $\gamma(-\infty)=q$ and $\gamma(0)\in S$.

Choosing for instance a Morse function on the manifold $S$ and adding a quadratic function in the normal direction this is easily seen to define
the \Poincare dual of the class $[S]\in\H_s(M;\Z/2)$ of $S$, i.e.~$\rho(S)=\PD[S]$.
The proof is immediately adjusted to the \Poincare dual case.
This concludes the proof of theorems \ref{thm:basic_theorem} and \ref{thm:basic_theorem_with_low_energy}.
\end{proof}
\section{Applications}
\label{sec:applications}

\subsection{Extrinsic topology of Lagrangian submanifolds}

\begin{Def}\label{def:displaceable}
We call a closed Lagrangian submanifold $L$ \textit{displaceable} in the closed symplectic manifold $(M,\om)$
if there exists a Hamiltonian diffeomorphism $\Phi_H\in\Ham(M,\om)$ such that $L\cap\Phi_H(L)=\emptyset$.
In particular, $\Phi_H$ displaces a sufficiently small neighborhood $\U(L)$ of $L$.

The \textit{displacement energy} of $L$ is $e(L):=\inf\{||H||\mid L\cap\Phi_H(L)=\emptyset\}$, i.e.~$L$ is displaceable iff $e(L)<\infty$.
\end{Def}

\begin{Thm}\label{thm:displaceable_implies_empty_moduli_spaces}
In the situation of theorem \ref{thm:basic_theorem} we denote by $\iota_k:\H_k(L;\Z/2)\pf\H_k(M;\Z/2)$
the homomorphism induced by the inclusion $\iota:L\subset M$.
If the Lagrangian submanifold $L$ is displaceable, then the homomorphism $\iota_k$ vanishes for degrees $k>\dim L+1-N_L$.
\end{Thm}
The homomorphism $\iota_0$ does not vanish, so that $N_L\leq\dim L+1$. This is well-known for displaceable Lagrangian submanifolds.
\begin{Cor}
For any displaceable monotone Lagrangian submanifold $L$ with $N_L\geq2$,
\beq
[L]=0\;\in\H_m(M;\Z/2)\,.
\eeq
\end{Cor}

\begin{proof}[\textsc{Proof of theorem \ref{thm:displaceable_implies_empty_moduli_spaces}}]
This is an application of theorem \ref{thm:basic_theorem}, namely we take a non-zero class $0\not=[S]\in\H_k(L;\Z/2)$ and represent its
image $\iota_k([S])\in\H_k(M;\Z/2)$ Floer theoretically with the help of theorem \ref{thm:basic_theorem}. Using the fact that $L$ is
displaceable we will prove that this cycle in Floer homology
vanishes simply by the fact that all moduli spaces involved are empty for a certain class of Hamiltonian functions. This shows that
$\PSS^{-1}\big(\iota_k([S])\big)=0$ and proves the assertion of theorem \ref{thm:displaceable_implies_empty_moduli_spaces}.

In contrast to theorem \ref{thm:basic_theorem} we denote here by $[S]$ a class in $\H_k(L;\Z/2)$ and by $\iota_k([S])$ its image in $\H_k(M;\Z/2)$.
We claim that under the assumption that $L$ is displaceable there exists a Hamiltonian function $G:S^1\times M\pf\R$ such that
$\M^-_{(L,S)}(\bar{x};J,G;\mathbf{0})=\emptyset$
for all $\bar{x}\in\Pt(G)$ of Conley-Zehnder index $-\CZ(\bar{x})<N_L$ and analogously for $\M^+_{(L,S)}(\bar{x};J,G;\mathbf{0})$.

By assumption $L$ is displaceable by the time-1-map $\Phi_K\in\Ham(M,\om)$ associated to a Hamiltonian function $K:S^1\times M\pf\R$, say.
In particular, $\Phi_K\big(\U(L)\big)\cap\U(L)=\emptyset$ for some small neighborhood $\U(L)$ of $L$.
We choose an autonomous Hamiltonian function
$H:M\pf\R_{\geq0}$ such that the support of $X_H$ is contained in $\U(L)$ and $H|_L>0$. Since $\Phi_K$ displaces the
support of $X_H$, we draw the standard conclusion
\beq
\Pt(\rho H\#K)=\Pt(K)
\eeq
for all constants $\rho\geq 0$ (cf.~\cite{Hofer_Zehner_Book}). Furthermore, $\A_{\rho H\#K}(\bar{x})=\A_K(\bar{x})$,
since by assumption $H\equiv0$ outside $\U(L)$ and all elements of $\Pt(K)$ and $\Pt(\rho H\#K)$ lie in $M\setminus\U(L)$.

We employ inequality \eqref{eqn:energy_estimate_M_(L,S)^-} to obtain the inequality
\beq
0\leq E(u)\;\leq\;-\A_{\rho H\#K}(\bar{x})-\inf_L{\big(\rho H\#K\big)}\;\leq\;-\A_K(\bar{x})-\rho\cdot\inf_L{H}-\inf_L{K}
\eeq
for an element $u\in\M^-_{(L,S)}(\bar{x};J,\rho H\#K;\mathbf{0})$.
We conclude that for $\rho\gg0$ the right hand side becomes negative, since $H|_L>0$. In particular,
all moduli spaces $\M^-_{(L,S)}(\bar{x};J,\rho H\#K;\mathbf{0})$ are empty for sufficiently large $\rho$.
This concludes the argument.
\end{proof}

In contrast to theorem \ref{thm:high_minimal_Maslov_number_implies_intersection} below it is essential that we assume that $L$ is displaced by a
Hamiltonian diffeomorphism and not a symplectic diffeomorphism as the examples $S^1\times\{\mathrm{pt}\}\subset S^1\times S^1$ shows.

Theorem \ref{thm:displaceable_implies_empty_moduli_spaces} was an application of theorem \ref{thm:basic_theorem}.
Using the same idea we present a new proof of Chekanov's result
\cite{Chekanov_Lagrangian_intersections_symplectic_energy_and_areas_of_holomorphic_curves} within our set-up
as an application of theorem \ref{thm:basic_theorem_with_low_energy}. We should mention that in
\cite{Chekanov_Lagrangian_intersections_symplectic_energy_and_areas_of_holomorphic_curves} this result is proved under the
sole assumption that $M$ is geometrically~bounded.

\begin{Thm}[Chekanov \cite{Chekanov_Lagrangian_intersections_symplectic_energy_and_areas_of_holomorphic_curves}]\label{thm:Chekanov}
Under the same assumptions as in theorem \ref{thm:basic_theorem}
the displacement energy $e(L)$ of a monotone closed Lagrangian submanifold $L$ with $N_L\geq2$ is bounded below by
the minimal area of a non-constant holomorphic disk or sphere: $e(L)\geq\min\{A_M,A_L\}$.
\end{Thm}

\begin{proof}
This goes along the same lines as the proof of the preceding theorem.
There is nothing to prove in case $e(L)=\infty$. Therefore, we assume from now on that $L$ is displaceable.

We choose a Hamiltonian function $K$ such that such that $\Phi_K^1$ displaces a small neighborhood $\U(L)$ of $L$.
Let us assume that the assertion of the theorem
is false, i.e.~$||K||<\min\big\{A_M,\,A_L\big\}$.

Again we choose a Hamiltonian function $H:M\pf\R_{\geq0}$ such that the support of $X_H$ is contained in $\U(L)$ and $H|_L>0$.
For an element $u\in\M^-_{(L,\mathrm{pt})}(\bar{x};J,\rho H\#K;\mathbf{0})$ we combine inequality \eqref{eqn:energy_estimate_M_(L,S)^-}
and the estimate of lemma \ref{Lemma_action_estimate_for_homologically_essential} for $\rho\geq0$ to obtain
\bea
E(u)&\;\leq\;-\A_{\rho H\#K}(\bar{x})-\inf_L(\rho H\#K)\\[0.5ex]
&\;\leq\;-\A_K(\bar{x})-\inf_L{K}-\rho\cdot\inf_L{H}\\[0.5ex]
&\;\leq\;\sup_MK-\inf_LK-\rho\cdot\inf_L{H}\\[0.5ex]
&\;\leq\;\sup_MK-\inf_MK\;=\;||K||\;\\[0.5ex]
&<\;\min\big\{A_M,\,A_L\big\}\,.
\eea
The second inequality is valid by exactly the same reasoning as in the proof of theorem \ref{thm:displaceable_implies_empty_moduli_spaces}.

Now we want to apply theorem \ref{thm:compactness}. Although the assumption $||\rho H\#K||<\min\big\{A_M,\,A_L\big\}$ does not hold
for large values of $\rho$ we established the crucial inequality $E(u)<\min\big\{A_M,\,A_L\big\}$ directly. As explained in remark
\ref{rmk:compactness_with_direct_inequality} we can still apply theorem \ref{thm:basic_theorem_with_low_energy} and thus
represent the class $[\mathrm{pt}]\in\H_0(M;\Z/2)$ in the Floer homology $\HF_*(\rho H\#K)$.

As in the proof of theorem \ref{thm:displaceable_implies_empty_moduli_spaces} we conclude from the displaceability of $L$ that
the moduli spaces $\M^-_{(L,\mathrm{pt})}(\bar{x};J,\rho H\#K;\mathbf{0})$ are empty if $\rho$ is sufficiently large.
This contradicts the fact that counting the number of elements of $\M^-_{(L,\mathrm{pt})}(\bar{x};J,\rho H\#K;\mathbf{0})$ defines
the class $[\mathrm{pt}]\not=0$.
Therefore, our assumption $||K||<\min\big\{A_M,\,A_L\big\}$ was false. This proves the theorem.
\end{proof}

\subsection{Lagrangian intersections and spectral capacities}$ $\\[1ex]
If we apply theorem \ref{thm:basic_theorem} in the special case $S=\mathrm{pt}$ to spectral capacities, we obtain.
\begin{Thm} \label{thm:high_minimal_Maslov_number_implies_intersection}
Let $L_0,L_1$ be two monotone closed Lagrangian submanifolds of the closed monotone symplectic manifold $(M^{2m},\om)$
of minimal Maslov number $N_{L_i}>\dim L_i+1$, $i=0,1$. If the spectral capacity $c_\gamma(M)$ of $M$ is finite,
then $L_0$ and $L_1$ intersect.
\end{Thm}
\begin{Rmk}
We note that theorem \ref{thm:high_minimal_Maslov_number_implies_intersection} is concerned with \textit{any} pair of Lagrangian
submanifolds. For example $L_1$ is not assumed to be a Hamiltonian deformation of $L_0$.
\end{Rmk}

\begin{proof}
As explained in remark \ref{rmk:extremal_cases_of_basic_theorems}, due to the assumption $N_{L_i}>\dim L_i+1$,
theorem \ref{thm:basic_theorem} provides an explicit
representation of $\PSS^{-1}([\om^m])$ in the Floer chain complex. We can use both, $L_0$ and $L_1$, for such a representation.

By definition, $c([\om^m],H)$ is the smallest action value of a representant of $[\om^m]$ in Floer homology. In particular,
from inequality \eqref{eqn:energy_estimate_M_(L,S)^-} we conclude
\beq
c([\om^m],H)\leq-\inf_{L_i}H\qquad\text{and}\qquad c(1,H)\geq-\sup_{L_i}H\,,\quad i=0,1
\eeq
using the property \textit{\Poincare duality} for the action selector in the second case. We obtain
\beq
\gamma(H)\geq\max\bigg\{\Big(-\sup_{L_1}H+\inf_{L_0}H\Big),\;\Big(-\sup_{L_0}H+\inf_{L_1}H\Big)\bigg\}\,.
\eeq
If the Lagrangian submanifolds do \textit{not} intersect, $L_0\cap L_1=\emptyset$, we can choose a sequence $\{H_n\}$ of Hamiltonian functions such that
$H_n|_{L_0}=n$ and $H_n|_{L_1}=0$.
In particular, $\gamma(H_n)\geq n$ and thus $c_\gamma(M)=\sup\{\gamma(H)\}=\infty.$ This proves the theorem.
\end{proof}

\begin{Cor}
All monotone closed Lagrangian submanifolds $L$ in a monotone closed symplectic manifold of finite
spectral capacity with minimal Maslov number $N_L>\dim L+1$ are connected. Furthermore,
$L\cap\varphi(L)\not=\emptyset$ $\forall\varphi\in\mathrm{Symp}(M,\om)$.
\end{Cor}

\begin{proof}
Each connected component of $L$ has minimal Maslov number at least $N_L$. Therefore, two such components of $L$ intersect.
For the second assertion is obvious from theorem \ref{thm:high_minimal_Maslov_number_implies_intersection}.
\end{proof}

\begin{Rmk}
Theorem \ref{thm:high_minimal_Maslov_number_implies_intersection} says that if the monotone symplectic $(M,\om)$ has two disjoint
Lagrangian submanifolds of sufficiently large minimal Maslov number, then it has infinite spectral capacity. 

The theorem allows for the following refinement: An arbitrarily small neighborhood of each such Lagrangian has infinite spectral capacity, namely
choose the sequence $\{H_n\}$ from the above proof such that $H_n$ is supported in the small neighborhood. In particular, $c_\gamma(L_i)=\infty$.
\end{Rmk}

\begin{Rmk}
\begin{enumerate}
\item The assumption of theorem \ref{thm:high_minimal_Maslov_number_implies_intersection} that the symplectic manifold has finite spectral
capacity is crucial as the example
$S^1\times\{\mathrm{pt}\}\subset S^1\times S^1$ shows: $S^1\times\{\mathrm{pt}\}$ is monotone and $N_{S^1\times\{\mathrm{pt}\}}=\infty$, but
$c_\gamma(S^1\times S^1)=\infty.$
\item The assumption on the minimal Maslov number is necessary, since otherwise the Lagrangian submanifold might be displaceable. Indeed,
there exist displaceable Lagrangian spheres $L$ in symplectic manifolds of the form $X^{n+1}\times\CP^n$,
cf.~\cite{Biran_Cieliebak_Symplectic_topology_on_subcritical_manifolds}.
Since $L$ is simply connected it is monotone and $N_L=\dim L+1$. Furthermore, an analogous calculation
as in lemma \ref{lemma:finite_spectral_capacity_of_CP^nxCP^n} shows $c_\gamma(X\times\CP^n)=c_\gamma(\CP^n)<\infty$ if $\om|_{\pi_2(X)}=0$.
\item We do not know whether the Lagrangian submanifold $L$ has to be monotone but we suspect that theorem
\ref{thm:high_minimal_Maslov_number_implies_intersection} does not generalize to non-monotone $L$. A counterexample could consist of an analog of two
$S^1\subset S^2$ where one of them is not an equator.
\end{enumerate}
\end{Rmk}

\begin{Thm}\label{thm:intersection_in_CP^nxCP^n}
Any two monotone Lagrangian submanifolds $L_0,L_1$ of $\CP^n\times\CP^n$ with minimal Maslov number $N_{L_i}> 2n+1\quad i=0,1$
have to intersect.
\end{Thm}

\begin{proof}
In view of theorem \ref{thm:high_minimal_Maslov_number_implies_intersection} we have to proof that $\CP^n\times\CP^n$ has finite spectral capacity.
This is content of lemma \ref{lemma:finite_spectral_capacity_of_CP^nxCP^n} and is proved with help of lemma
\ref{lemma:finite_spectral_capacity_of_CP^n}.
\end{proof}
\begin{Lemma}\label{lemma:finite_spectral_capacity_of_CP^n}
The monotone symplectic manifold $(\CP^n,\om_{\mathrm{FS}})$ has finite spectral capacity
\beq
\tilde{c}_{\tilde{\gamma}}(\CP^n)=c_\gamma(\CP^n)\leq\om_{\mathrm{FS}}(\CP^1)\,.
\eeq
\end{Lemma}
\begin{proof}
The equality is due to the fact that the minimal Chern number $N_{\CP^n}$ is sufficiently large,
namely $2N_{\CP^n}>\dim\CP^n$. We prove only the inequality.

The quantum cohomology ring of $\CP^n$ over $\Novo=\Z/2$ is given by
\beq
\QH^*(\CP^n)\cong\frac{(\Z/2)[\!\!\:[q]\!\!\:][q^{-1}][p]}{<p^{n+1}=q>}\,,
\eeq
i.e. Laurent series in $q$ and polynomials in $p$ up to degree $n$, where the class of the symplectic form
corresponds to $p$. Furthermore, $q$ corresponds to the class $[\CP^1]$ of the holomorphic sphere $\CP^1$. In particular,
the relation $p^{n+1}=q$ reads
\beq
[\om_{\mathrm{FS}}]\ast[\om_{\mathrm{FS}}^n]=1_{\CP^n}\ast e^{[\CP^1]}\,.
\eeq
See \cite{McDuff_Salamon_J_holomorphic_curves_and_symplectic_topology} for details on the quantum cohomology ring of $\CP^n$.

Using the properties \textit{Novikov} and \textit{Product} of the action selector we find
\bea
c([\om_{\mathrm{FS}}]\ast[\om_{\mathrm{FS}}^n], H)&\;=\;c(1_{\CP^n}\ast e^{[\CP^1]},H)=c(1_{\CP^n},H)-\om_{\mathrm{FS}}(\CP^1)\quad\text{and}\\[1ex]
c([\om_{\mathrm{FS}}]\ast[\om_{\mathrm{FS}}^n], H)&\;\leq\;c([\om_{\mathrm{FS}}],0)+c([\om_{\mathrm{FS}}^n], H)=c([\om_{\mathrm{FS}}^n], H)\,.
\eea
This implies $\gamma(H)=c(1_{\CP^n},H)-c([\om_{\mathrm{FS}}^n], H)\leq\om_{\mathrm{FS}}(\CP^1)$,
and thus $c_\gamma(\CP^n)\leq\om_{\mathrm{FS}}(\CP^1)$.
\end{proof}

\begin{Lemma}\label{lemma:finite_spectral_capacity_of_CP^nxCP^n}
The spectral capacity of $(\CP^n\times\CP^n,\om_{\mathrm{FS}}\oplus\om_{\mathrm{FS}})$ is finite:
\beq
c_\gamma(\CP^n\times\CP^n)\leq2\om_{\mathrm{FS}}(\CP^1)\,.
\eeq
\end{Lemma}
\begin{proof}
We abbreviate
\beq
a:=[\om_{\mathrm{FS}}]\oplus0\,,\zw b:=0\oplus[\om_{\mathrm{FS}}]\,,\zw
A:=\CP^1\times\{\mathrm{pt}\}\zw\text{and}\zw B:=\{\mathrm{pt}\}\times\CP^1\,.
\eeq
The class $a^nb^n:=(a\cup\cdots\cup a)\cup(b\cup\cdots\cup b)$ represents (up to a factor $\binom{2n}{n}$) the class
$[(\om_{\mathrm{FS}}\oplus\om_{\mathrm{FS}})^{2n}]$.
Since both factors in $\CP^n\times\CP^n$ have the same symplectic form we can use the Kuenneth formula for quantum cohomology
\cite[chapter 11.1]{McDuff_Salamon_J_holomorphic_curves_and_symplectic_topology}.
Therefore, together with lemma \ref{lemma:finite_spectral_capacity_of_CP^n} we compute in the quantum cohomology of $\CP^n\times\CP^n$:
\beqn
(ab)\ast(a^nb^n)\;=\;(a\ast a^n)(b\ast b^n)\;=\;(1_{\CP^n}\ast e^A)(1_{\CP^n}\ast e^B)\;=\;1_{\CP^n\times\CP^n}\ast e^{(A+B)}\,.
\eeq
Now the reasoning is as above, namely
\bea
c\big((ab)\ast(a^nb^n),H\big)&\;=\;c(1_{\CP^n\times\CP^n}\ast e^{(A+B)},H)\;=\;c(1_{\CP^n\times\CP^n},H)-2\om_{\mathrm{FS}}(\CP^1)\\[1ex]
c\big((ab)\ast(a^nb^n),H\big)&\;\leq\;c(ab,0)+c(a^nb^n,H)\;=\;c(a^nb^n,H)\,.
\eea
Note, that $c(\binom{2n}{n}\cdot a^nb^n,H)=c(a^nb^n,H)$. The lemma follows.
\end{proof}

\begin{Cor}\label{cor:simply_connected_in_CP^nxCP^n_intersect}
Any two simply connected Lagrangian submanifolds in $\CP^n\times\CP^n$ intersect each other.
Furthermore, they are all connected.
\end{Cor}

\begin{proof}
Since a simply connected Lagrangian submanifold $L$ in a monotone symplectic manifold is monotone with minimal Maslov number equal to twice the minimal
Chern number, we conclude $N_L=2N_{\CP^n\times\CP^n}=2N_{\CP^n}=2(n+1)>2n+1$ and can apply theorem \ref{thm:intersection_in_CP^nxCP^n}.
\end{proof}

\begin{Ex}
An example of a simply connected Lagrangian submanifold of $\CP^n\times\CP^n$ is the anti-diagonal $\bar{\Delta}=\{(\bar{z},z)\mid z\in\CP^n\}$.
Unfortunately, we are not aware of any other examples of simply connected Lagrangian submanifolds of $\CP^n\times\CP^n$ or of non-simply connected
monotone Lagrangian submanifolds of sufficiently large minimal Maslov number.
\end{Ex}

\begin{Rmk}
Biran informed us that using his techniques he can prove corollary \ref{cor:simply_connected_in_CP^nxCP^n_intersect}, too.
Using tools adapted to four dimensions, Hind \cite{Hind_Lagrangian_spheres_in_S^2xS^2} proves that all Lagrangian spheres
in $S^2\times S^2$ are Hamiltonianly isotopic to the anti-diagonal. Since the minimal Maslov number of such spheres equals 4 they are not displaceable
and corollary \ref{cor:simply_connected_in_CP^nxCP^n_intersect} follows in the case $n=1$.
\end{Rmk}

Cornea \cite{Cornea_slides_SMS_2004} reported on an ongoing project with Lalonde at the ''SMS 2004''
in Montreal\footnote{see http://www.dms.umontreal.ca/sms/index.html}. The following theorem is a special case of their results. We give here a
new proof since it demonstrates our method quite nicely. Another proof of theorem \ref{thm:Cornea_Lalonde}, in spirit of Gromov's theorem
asserting that $\H^1(L;\R)\not=0$ for all closed Lagrangian $L\subset\R^{2n}$, was explained to us by Salamon.

\begin{Thm}\label{thm:Cornea_Lalonde}
Let $(M,\om)$ be a symplectically aspherical manifold and $L\subset M$ a monotone closed Lagrangian submanifold with $N_L\geq2$.
If $L$ is displaceable, then through \emph{each} point of $L$ passes a non-constant holomorphic disk
with boundary lying entirely on $L$ and which is of Maslov index $\Mas\leq\dim L+1$.
\end{Thm}

\begin{proof}
We argue by contradiction and assume that there is a point $p_0\in L$ through which passes \textit{no} holomorphic disk of Maslov index
$\Mas\leq\dim L+1$. Then inspection of the proof of theorem \ref{thm:compactness} shows that the moduli spaces
$\M^\pm_{(L,\{p_0\})}(\bar{x};J,H;\mathbf{0})$ are compact up to splitting of multiple Floer cylinders for \textit{all} Hamiltonian functions $H$.
In other words, no bubbling-off of holomorphic disks can occur. Furthermore, no holomorphic spheres exist since $M$ is assumed to be symplectically
aspherical.

Now we argue as in the proof of theorem \ref{thm:displaceable_implies_empty_moduli_spaces} or \ref{thm:Chekanov}.
Namely, since $L$ is displaceable, we can show that the moduli spaces $\M^\pm_{(L,\{p_0\})}(\bar{x};J,H;\mathbf{0})$ are empty for appropriately
chosen Hamiltonian function $H$. On the other hand, we can represent the non-vanishing class $[p_0]\in\H_0(M;\Z/2)$ by
counting the number of solutions. This contradiction finishes the proof.
\end{proof}

Since the dimension of $\M^\pm_{(L,\{p_0\})}(\bar{x};J,H;\mathbf{0})$ equals $\pm\CZ(\bar{x})-\dim L$ we can relax the assumption
of $(M,\om)$ being symplectically aspherical to $N_M>\dim L$. In this case there might be holomorphic spheres but the index formula
\eqref{eqn:Fredholm_index_formula_bubbling} shows that the moduli spaces $\M^\pm_{(L,\{p_0\})}(\bar{x};J,H;\mathbf{0})$ are compact
for all relevant $\bar{x}$, more precisely for $\CZ(\bar{x})=\pm\dim M$. For instance the theorem holds for $\CP^n$.

Combining the proofs of theorem \ref{thm:high_minimal_Maslov_number_implies_intersection} and \ref{thm:Cornea_Lalonde} we obtain the following assertion
for a symplectic manifold $(M,\om)$ of finite spectral capacity.

If there exists a point on the monotone Lagrangian submanifold $L\subset M$ through which passes \textbf{no} non-constant holomorphic disk of Maslov
index $\Mas\leq\dim L+1$, e.g.~$N_L>\dim L+1$, then $L\cap\varphi(L)\not=\emptyset$ for all $\varphi\in\mathrm{Symp}(M,\om)$.

The example $S^1\subset S^2$ shows that this statement can't be reversed. Through each point of an equator passes a holomorphic disk of Maslov index 2
and each image under a symplectomorphism intersect the equator again. More generally, this holds true for the Clifford torus in $\CP^n$, see
\cite{Entov_Polterovich_Quasi-states_and_symplectic_intersections}.
%
\bibliographystyle{amsplain}
\bibliography{../Bibtex/bibtex_paper_list}
\end{document}